\documentclass{article}[12pt]  
\usepackage[utf8x]{inputenc}
\usepackage{amsfonts} 
\usepackage{graphics}
\usepackage{graphicx} 
\usepackage{color}
\usepackage[margin=3.cm]{geometry}

\usepackage{amsmath,amssymb}
\usepackage{color}
\usepackage[dvipsnames]{xcolor}
\usepackage{appendix}
\usepackage{amsthm}
\usepackage{array}
\usepackage{mathtools}

\renewcommand{\epsilon}{\varepsilon}

\graphicspath{{figures/}}

\DeclareGraphicsExtensions{.eps,.pdf,.png,.jpg,.JPG,.jpeg,.ps,.pdf_tex}
\numberwithin{equation}{section}
\begin{document}


  \title{Asymptotic analysis for close evaluation of layer potentials}
\author{Camille Carvalho\footnote{Applied Mathematics Unit, School of Natural Sciences,
    University of California, Merced, 5200 North Lake Road, Merced, CA
    95343}  \and  Shilpa Khatri ${}^{\ast}$  \and  Arnold D. Kim ${}^{\ast}$}
   \maketitle

  \begin{abstract}
    We study the evaluation of layer potentials close
      to the domain boundary. Accurate evaluation of layer potentials near boundaries
    is needed in many applications, including fluid-structure
    interactions and near-field scattering in nano-optics. 
    When numerically evaluating layer potentials, it is natural to use
    the same quadrature rule as the one used in the Nystr\"om method
    to solve the underlying boundary integral equation. However, this
    method is problematic for evaluation points close to
    boundaries. For a fixed number of quadrature points, $N$, this
    method incurs $O(1)$ errors in a boundary layer of thickness
    $O(1/N)$. Using an asymptotic expansion for the kernel of the
    layer potential, we remove this $O(1)$ error. We demonstrate the
    effectiveness of this method for interior and exterior problems
    for Laplace's equation in two dimensions.
  \end{abstract}



\section{Introduction}
\label{sec:intro}

Boundary integral equation methods are useful for solving boundary
value problems for linear, elliptic, partial differential equations
~\cite{guenther1996partial, kress1999linear}. Rather than solving the
partial differential equation directly, one represents the solution as
a layer potential, an integral operator applied to a density. The
density is the solution of an integral equation on the boundary of the
domain that includes the prescribed boundary data.  This formulation
offers several advantages for the numerical solution of boundary value
problems. First, the dominant computational cost is from the integral
equation on the boundary whose dimension is lower than that of the
domain. Second, this boundary integral equation can be solved to very
high order using Nystr\"om methods~\cite{atkinson1997secondkind,
  delves1988computational}. Finally, the solution, given by this layer
potential, can be evaluated anywhere in the domain without restriction
to a particular mesh. For these reasons, boundary integral equations
have found broad applications, including in fluid mechanics and
electromagnetics.

One particular challenge in using boundary integral equation methods
is the so-called close evaluation problem~\cite{barnett2014evaluation,
  helsing2008evaluation}.  Since a layer potential is an integral over
the boundary, it is natural to evaluate it numerically using the same
quadrature rule used in the Nystr\"om method to solve the boundary
integral equation. In that case, we say that the layer potential is
evaluated using its native quadrature rule. Numerical evaluation of
the layer potential using its native quadrature rule inherits the high
order accuracy associated with solving the boundary integral equation,
except for points close to the boundary. For these close evaluation
points, the native quadrature rule produces an $O(1)$ error.  This
$O(1)$ error is due to the sharply peaked kernel of
  the layer potential leading to its nearly singular behavior.

%
%
%
%
There are several problems that require accurate layer
  potential evaluations close to the boundary of the domain. For
example, modeling of micro-organisms swimming, suspensions of
droplets, and blood cells in Stokes flow use boundary integral methods
~\cite{smith2009boundary, barnett2015spectrally, marple2016fast,
  keaveny2011applying}. The key to these problems is the accurate
computation of velocity fields or forces close to the boundary as
these quantities provide the physical mechanisms leading to locomotion
and other phenomena of interest.  Another example is in the field of
plasmonics~\cite{Maier07}, where one seeks to gain control of light at
the sub-wavelength scales for applications such as
nano-antennas~\cite{akselrod2014probing,novotny2011antennas} and
sensors~\cite{mayer2008label,sannomiya2008situ}.  Surface plasmons are
sub-wavelength fields localized at interfaces between nano-scale metal
obstacles and their surrounding dielectric background medium.  Thus,
these problems require accurate computation of electromagnetic fields
near interfaces. These problems and others motivate
  the need to address the close evaluation problem.

The close evaluation problem for layer potentials has been studied
previously for Laplace's equation.  For example, Beale and
Lai~\cite{beale2001method} have studied this problem in two dimensions
by first regularizing the nearly singular kernel and then adding
corrections for both the discretization and the regularization. The
result of this approach is a uniform error in space.  This method has
been extended to three-dimensional problems~\cite{beal2016asimple}.
Helsing and Ojala~\cite{helsing2008evaluation} have studied the
Dirichlet problem in two dimensions by combining a globally
compensated quadrature rule along with interpolation to achieve very
accurate results over all regions of the domain.
Barnett~\cite{barnett2014evaluation} has also studied this problem in
two dimensions. In that work, Barnett has established a bound for the
error associated with the periodic trapezoid rule.  We make use of
this result in our work below.  To address the $O(1)$ error in the
close evaluation problem, Barnett has used surrogate local expansions
with centers placed near, but not on, the boundary. 
  This new method, called quadrature by expansion (QBX), leads to very
  accurate evaluations of the layer potential close to the
  boundary. Further error analysis of this method and extensions to
  evaluations on the boundary for the Helmholtz equation is presented in
  Kl\"{o}ckner {\it et al}~\cite{klockneretal2013}.  
 For the special case of rectangular domains, Fikioris
  {\it et al}~\cite{fikioris1987strongly, fikioris1988strongly} have
  addressed the close evaluation problem by removing problematic terms
  from the explicit eigenfunction expansion of the fundamental
  solution.  

Here, we develop a new method to address the close
  evaluation problem. We first determine the asymptotic behavior of
  the sharply peaked kernel and then use that to approximate the layer
  potential. Doing so relieves the quadrature rule from having to
  integrate over a sharp peak. Instead, the quadrature rule is used to
  correct the error made by this approximation.  This new method is
  accurate, efficient, and easy to implement.

In this paper, we study the close evaluation problem in two dimensions
for Laplace's equation. We use a Nystr\"om method based on the
periodic trapezoid rule to solve the boundary integral equation.  We
study the double-layer potential for the interior Dirichlet problem
and the single-layer potential for the exterior Neumann problem. For
both of these problems, we introduce an asymptotic expansion for the
sharply peaked kernel of the layer potential for close evaluation
points, which is the main cause for error.  Using the Fourier series
of this asymptotic expansion, we compute its contribution to the layer
potential with spectral accuracy.  Through several examples, we show
that this asymptotic method effectively reduces errors in the close
evaluation of layer potentials.

The remainder of this paper is as follows. In Section \ref{sec:circle}
we study in detail the illustrative example of the interior Dirichlet
problem for a circular disk. For this problem, we obtain an explicit
error when using the periodic trapezoid rule to evaluate the
double-layer potential. This error motivates the use of an asymptotic
expansion for the sharply peaked kernel in the general method we
develop in Section \ref{sec:doublelayer} to evaluate the double-layer
potential for the interior Dirichlet problem. In Section
\ref{sec:singlelayer}, we extend this method to the single-layer
potential for the exterior Neumann problem. We discuss the general
implementation of these methods in Section
\ref{sec:implementation}. Section
  \ref{sec:conclusions} gives our conclusions.

\section{Illustrative example: interior Dirichlet problem for a
  circular disk}
\label{sec:circle}

We first study the close evaluation problem for 
\begin{subequations}
  \begin{gather}
    \Delta u = 0 \quad \text{in $D = \{ r < a \}$}, \label{eq:2.1a}\\
    u=f \quad \text{on $\partial D = \{ r = a \}$}. \label{eq:2.1b}
  \end{gather}
  \label{eq:2.1}
\end{subequations}
For this problem, we compute an explicit error when
  applying an $N$-point periodic trapezoid rule (PTR$_{N}$). This
  error reveals the key factors leading to the large errors observed
  in the close evaluation problem. Moreover, this analysis provides
  the motivation for the general asymptotic method.

It is well understood that the solution of \eqref{eq:2.1} is given by
Poisson's formula~\cite{strauss1992partial}. Here, we seek the
solution as the double-layer potential~\cite{kress1999linear},
\begin{equation}
  u(\mathbf{x}) = \frac{1}{2\pi}  \int_{|\mathbf{y}| = a}
  \frac{\mathbf{n}_{y} \cdot (\mathbf{x} - \mathbf{y})}{|\mathbf{x} -
    \mathbf{y}|^{2}} \mu(\mathbf{y}) \mathrm{d}\sigma_{y}.
  \label{eq:2.2}
\end{equation} 
Here, $\mathbf{x} \in D$ denotes the evaluation point,
$\mathbf{y} \in \partial D$ denotes the variable of integration,
$\mathbf{n}_{y}$ denotes the unit, outward normal at $\mathbf{y}$, and
$\mathrm{d}\sigma_{y}$ denotes a differential boundary element. The
density, $\mu(\mathbf{y})$, satisfies the following boundary integral
equation,
\begin{equation}
  - \frac{1}{2} \mu(\mathbf{y}) - \frac{1}{4\pi a} \int_{|\mathbf{y}'
    | = a} \mu(\mathbf{y}') \mathrm{d}\sigma_{y'} = f(\mathbf{y}),
  \label{eq:2.3}
\end{equation}
from which we determine that
\begin{equation}
  \mu(\mathbf{y}) = \frac{1}{2\pi a} \int_{| \mathbf{y}' | = a}
  f(\mathbf{y}') \mathrm{d}\sigma_{y'} - 2 f(\mathbf{y}).
  \label{eq:2.4}
\end{equation}

To numerically evaluate \eqref{eq:2.2}, we substitute the
parameterization, $\mathbf{x} = (r \cos t^{\ast}, r \sin t^{\ast})$,
and $\mathbf{y} = (a \cos t, a \sin t)$, with $r < a$, and
$t^{\ast}, t \in [0,2\pi]$, and obtain
\begin{equation}
  u(r,t^{\ast}) = \frac{1}{2\pi}  \int_0^{2\pi}  
  \left[ \frac{ a r  \cos ( t - t^{\ast} ) - a^{2}}{a ^{2}+ r^2 - 2
      ar \cos( t - t^{\ast} ) } \right] \mu(t) \mathrm{d}t =
  \frac{1}{2\pi} \int_{0}^{2\pi} K(t - t^{\ast}) \mu(t) \mathrm{d}t.
  \label{eq:2.5}
\end{equation} 
Using PTR$_{N}$ with points $t_{j} = 2\pi j/N$ for
$j = 1, \cdots, N$, to approximate \eqref{eq:2.5}, we obtain
\begin{equation}
  u(r,t^{\ast}) \approx U^{N}(r,t^{\ast}) = \frac{1}{N} \sum_{j =
    1}^{N} K\left( \frac{2\pi j}{N} - t^{\ast} \right) \mu\left(
    \frac{2\pi j}{N} \right).
  \label{eq:2.6}
\end{equation}

\begin{figure}[t]
  \centering
  \includegraphics[width=0.4\linewidth]{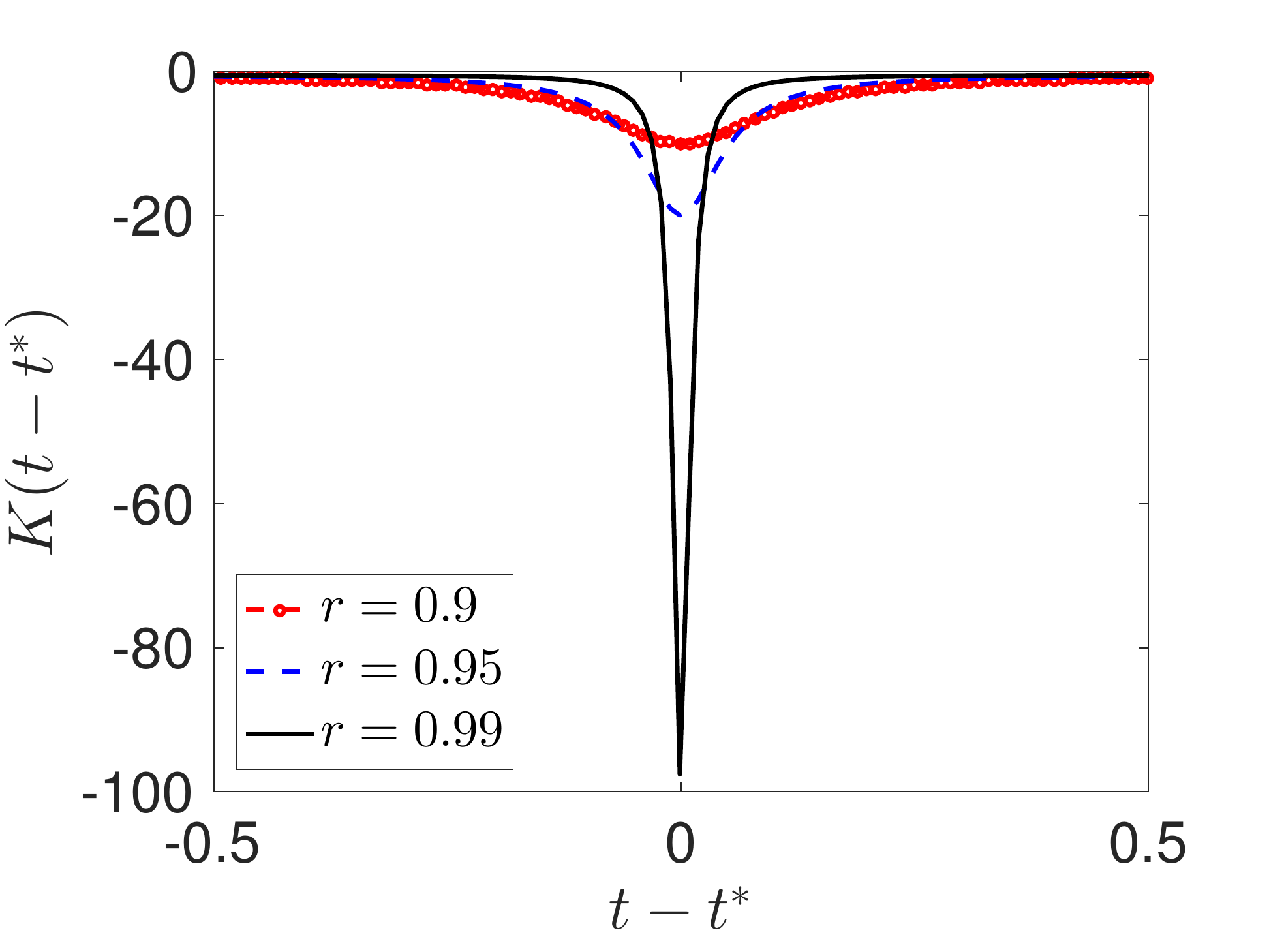}
  \includegraphics[width=0.4\linewidth]{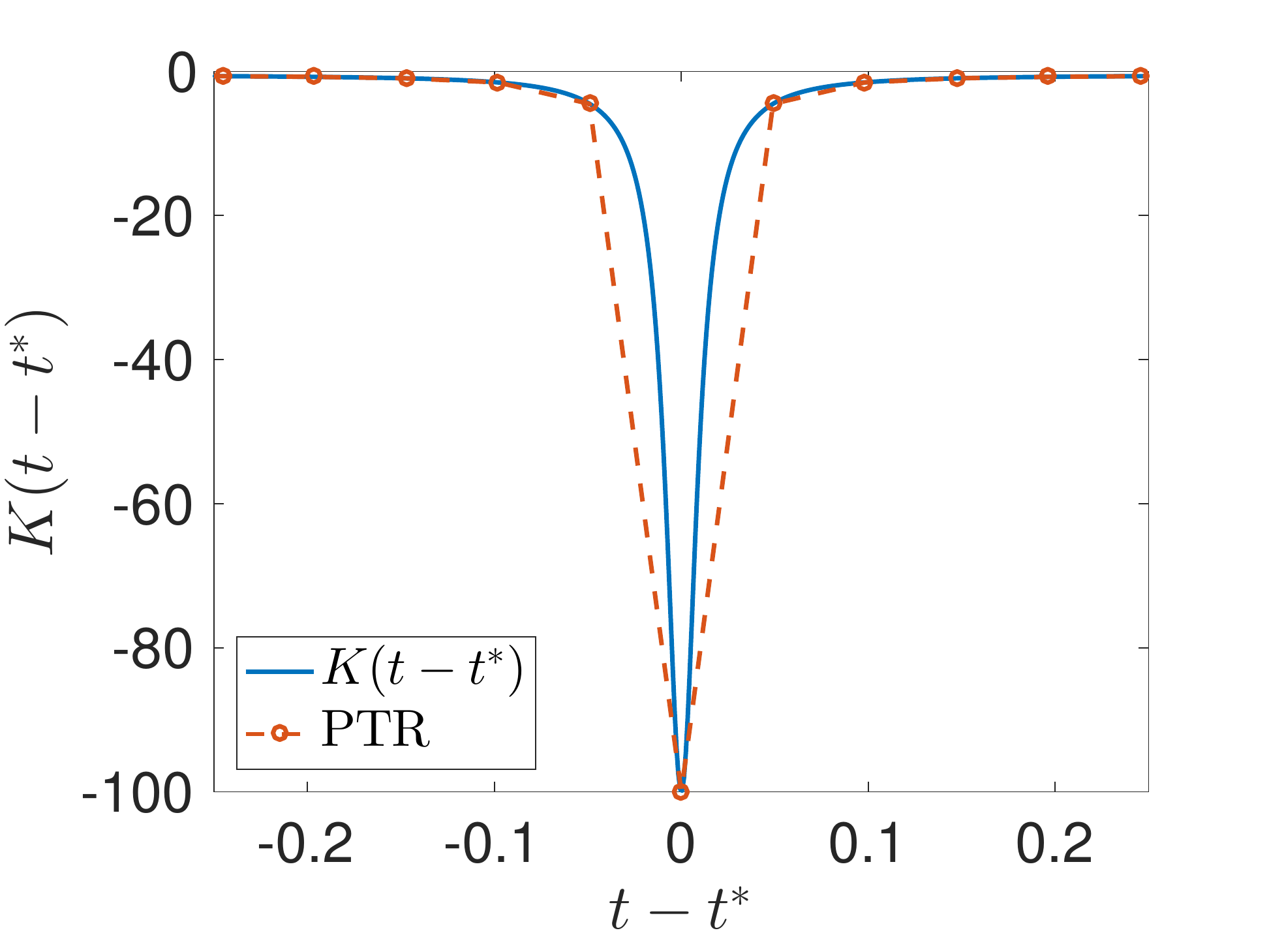}
  \caption{[Left] Plot of the kernel, $K(t - t^{\ast})$, defined in
    \eqref{eq:2.5} with $a = 1$ for $r = 0.9$ (dot-dashed curve),
    $0.95$ (dashed curve), and $0.99$ (solid curve). [Right] Plot of
    the kernel, $K(t - t^{\ast})$, with $a = 1$, and $r = 0.99$ (solid
    curve), and the corresponding piecewise linear approximation
    associated with PTR$_{128}$ (dot-dashed curve).}
  \label{fig:1}
\end{figure}

The error, $U^{N}(r,t^{\ast}) - u(r,t^{\ast})$, is not uniform in $D$.
In particular, $U^{N}$ is very accurate for evaluation points far away
from the boundary. On the other hand, it is $O(1)$ when $r \sim a$.
The reason for this large error is due to the kernel,
$K(t - t^{\ast})$.  In Fig.~\ref{fig:1}, we show plots of $K$ as a
function of $t - t^{\ast}$. The left plot of Fig.~\ref{fig:1} shows
that $K$ becomes sharply peaked about $t = t^{\ast}$ when $r \sim a$.
We do not evaluate the double-layer potential on the
boundary. Nonetheless, because $K$ becomes sharply peaked as
$r \to a$, we say that the double-layer potential is nearly singular.
  The right plot of Fig.~\ref{fig:1} shows that the piecewise linear
  approximation of $K$ used by PTR$_{128}$ will grossly overestimate
  the magnitude of the double-layer potential for $r = 0.99 a$.  It is
  this error that leads to the $O(1)$ errors produced by PTR$_{N}$ for
  close evaluation points.  Barnett~\cite{barnett2014evaluation} has
  shown that this error is $O(1)$ for $a - r = O(1/N)$. In light of
  this result, we say that the error exhibits a boundary layer of
  thickness $O(1/N)$ in which it undergoes rapid growth.

  In the limit as $N \to \infty$, PTR$_{N}$ converges because the
  boundary layer vanishes at a rate of $O(1/N)$. However, that is not
  the limit we consider here. Rather, we study the limit as the
  evaluation point approaches the boundary with $N$ fixed. For that
  case, PTR$_{N}$ is unable to accurately capture the sharp peak of
  the kernel about $t = t^{\ast}$ that forms as $r \to a$. 

Using the error associated with
PTR$_{N}$~\cite{davis1959ptr}, we find $U^{N}$
defined in \eqref{eq:2.6} satisfies
\begin{equation}
  U^{N}(r,t^{\ast}) = u(r,t^{\ast}) + \sum_{\substack{l = -\infty\\l
      \neq 0}}^{\infty} \hat{p}[lN],
  \label{eq:2.7}
\end{equation}
where
\begin{equation}
  \hat{p}[k] = \frac{1}{2\pi} \int_{0}^{2\pi} K(t - t^{\ast}) \mu(t)
  e^{-\mathrm{i} k t } \mathrm{d}t.
  \label{eq:2.8}
\end{equation}
The error in \eqref{eq:2.7} is aliasing of high frequencies.  To
determine $\hat{p}[k]$ explicitly, we use the Fourier series
representation of the kernel,
\begin{equation}
  K(t-t^{\ast}) = \frac{ ar  \cos ( t - t^{\ast} )
    -a^2 }{a ^{2}+ r^2 - 2  ar\cos ( t - t^{\ast} )} =
  -\frac{1}{2} -  \frac{1}{2} \sum_{m = -\infty}^{\infty}
  \left( \frac{r}{a} \right)^{|m|} e^{\mathrm{i} m (t - t^{\ast})},
  \label{eq:2.9}
\end{equation} 
and of the density
\begin{equation}
  \mu(t) = \sum_{n = -\infty}^{\infty} \hat{\mu}[n] e^{\mathrm{i} n t},
  \label{eq:2.10}
\end{equation}
to find 
\begin{equation}
  K(t - t^{\ast}) \mu(t) = -\frac{1}{2} \sum_{n =
    -\infty}^{\infty} \hat{  {\mu}}[n] e^{\mathrm{i} n t} -\frac{1}{2}
  \sum_{m = -\infty}^{\infty} \left( \frac{r}{a} \right)^{|m|}
  e^{-\mathrm{i} m t^{\ast}} \sum_{n = -\infty}^{\infty} \hat{\mu}[n]
  e^{\mathrm{i} (m + n) t}.
  \label{eq:2.11}
\end{equation}
Substituting \eqref{eq:2.11} into \eqref{eq:2.8}, and
rearranging terms, we find that
\begin{equation}
  \hat{p}[k] = -\hat{\mu}[k] -\frac{1}{2} \sum_{m = 1}^{\infty} \left(
    \frac{r}{a} \right)^{m} \left( e^{\mathrm{i} m t^{\ast}}
    \hat{\mu}[k+m] + e^{-\mathrm{i} m t^{\ast}} \hat{\mu}[k-m] \right).
  \label{eq:2.12}
\end{equation}

\begin{figure}[t]
  \centering
 \includegraphics[width=0.7\linewidth]{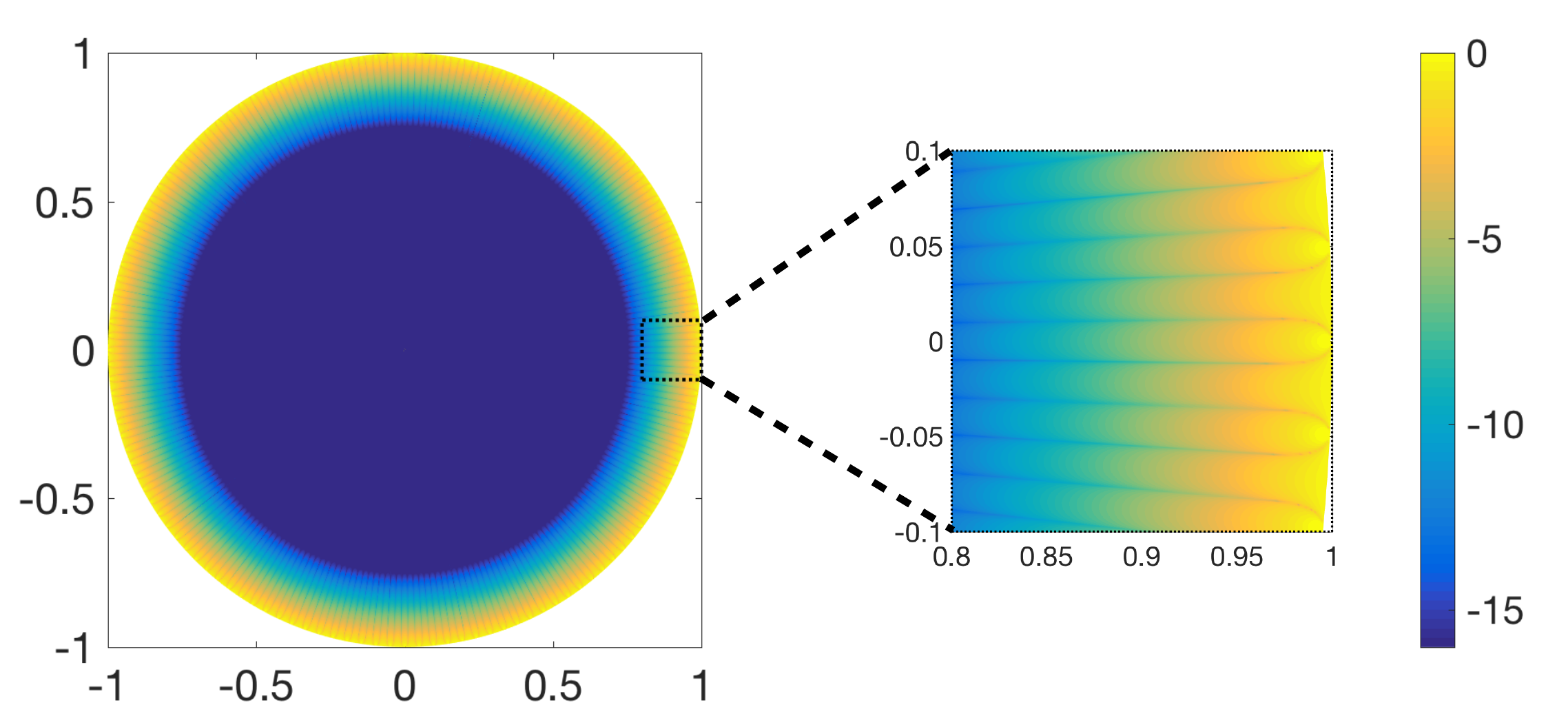}
  
 \caption{Contour plot of $\log_{10} |E^{N}(r,t)|$ where $E^{N}$ is
   given by \eqref{eq:2.15} with $a = 1$ and $N = 128$.}
  \label{fig:2}
\end{figure}

We now obtain the error in using PTR$_{N}$ to
evaluate the double-layer potential by substituting \eqref{eq:2.12}
into \eqref{eq:2.7} which yields
\begin{multline}
  E^{N}(r,t^{\ast}) = U^{N}(r,t^{\ast}) - u(r,t^{\ast}) = \sum_{l =
    1}^{\infty} \left\{ -\hat{\mu}^{\ast}[lN] - \frac{1}{2} \sum_{m =
      1}^{\infty} \left[ \left( \frac{r}{a} \right)^{m} \left(
        \hat{\mu}^{\ast}[lN-m] e^{\mathrm{i} m t^{\ast}} +
        \hat{\mu}^{\ast}[lN+m]
        e^{-\mathrm{i} m t^{\ast}} \right) \right] \right\}\\
  + \sum_{l = 1}^{\infty} \left\{ - \hat{\mu}[lN] -\frac{1}{2} \sum_{m
      = 1}^{\infty} \left[ \left( \frac{r}{a} \right)^{m} \left(
        \hat{\mu}[lN+m] e^{\mathrm{i} m t^{\ast}} + \hat{\mu}[lN-m]
        e^{-\mathrm{i} m t^{\ast}} \right) \right] \right\}.
  \label{eq:2.13}
\end{multline}
Here, we have assumed $\mu$ is real, so
$\hat{\mu}[-k] = \hat{\mu}^{\ast}[k]$, where $[ \cdot ]^{\ast}$
denotes complex conjugation. Suppose we have chosen $N$ to be large
enough so that $\hat{\mu}[lN] \ll 1$ for $l > 0$.
For that case, only terms in \eqref{eq:2.13} proportional to
$\hat{\mu}[lN-m]$ will substantially contribute to the error. By
neglecting the other terms, we obtain
\begin{equation}
  E^{N}(r,t^{\ast}) \sim - \sum_{l = 1}^{\infty} \sum_{m
    = 1}^{\infty} \left( \frac{r}{a} \right)^{m} \left[ \text{Re}\{
    \hat{\mu}[lN-m] \} \cos( m t^{\ast} ) + \text{Im}\{
    \hat{\mu}[lN-m] \} \sin( m t^{\ast} ) \right].
  \label{eq:2.14}
\end{equation}

Equation \eqref{eq:2.14} is the asymptotic error made by
PTR$_{N}$. The key point is that
  when $N$ is fixed, 
this error is not uniform for $r \in [0,a)$. When $r \ll a$, we see
that the error is much smaller than when $r \sim a$.  Consider the
specific case in which $\mu = 1$, so that $\hat{\mu}[0] = 1$ and
$\hat{\mu}[k] = 0$ for all $k \neq 0$. For that case, \eqref{eq:2.14}
simplifies to
\begin{equation}
  E^{N}(r,t^{\ast}) \sim \frac{\left( \frac{r}{a} \right)^{2N} -
    \left( \frac{r}{a} \right)^{N} \cos (N t^{\ast})}{1 + \left(
      \frac{r}{a} \right)^{2N} - 2 \left( \frac{r}{a} \right)^{N} \cos
    (N t^{\ast})}. 
  \label{eq:2.15}
\end{equation}
According to \eqref{eq:2.15},
$| E^{N}(a(1-\epsilon),t^{\ast}) | = O((1-\epsilon)^{N}) =
O(e^{-\epsilon N})$,
and $| E^{N}(r,t^{\ast}) | \to 1/2$ as $r \to a$.  These results show
that $E^{N}$ has a boundary layer of thickness $O(1/N)$ where it
exponentially increases to values that are $O(1)$. Fig.~\ref{fig:2}
shows a plot of \eqref{eq:2.15} over the entire circular disk and a
close-up near the boundary.  These plots show the boundary layer about
$r = a$ where the error attains $O(1)$ values.  In practice, we would
like to set $N$ based on the resolution required to solve the boundary
integral equation. It is neither desirable nor practical to increase
$N$ just to reduce aliasing in the evaluation of the double-layer
potential.  In light of this, we make the following observations.
\begin{itemize}

\item Equation \eqref{eq:2.13} gives the error incurred by
PTR$_{N}$ to approximate the double-layer
  potential. This error is due to aliasing. Equation \eqref{eq:2.14}
  gives the asymptotic approximation of this error when the $N$-point
  grid sufficiently samples $\mu$.

\item The aliasing error is not uniform with respect to $r$. For the
  case in which $\mu = 1$, the asymptotic error simplifies to
  \eqref{eq:2.15}. From this result, we find that the error grows
  rapidly and becomes $O(1)$ in a boundary layer of thickness $O(1/N)$
  near the boundary. This boundary layer is shown in Fig.~\ref{fig:2}.

\item For points within the boundary layer, the sharply peaked kernel
  causes aliasing due to insufficient resolution.  Fig.~\ref{fig:1}
  shows how the sharp peak of the kernel when $r/a = 0.99$ is
  under-resolved on the grid for PTR$_{128}$.

\end{itemize}

Alternatively, by substituting \eqref{eq:2.9} and
  \eqref{eq:2.10} into \eqref{eq:2.5}, we obtain
\begin{equation}
  u(r,t^{\ast}) = \sum_{n = -\infty}^{\infty} \hat{K}^{\ast}[n]
  \hat{\mu}[n] e^{- \mathrm{i} n t^\ast} \approx \sum_{n =
    -N/2}^{N/2-1} \hat{K}^{\ast}[n] \hat{\mu}[n] e^{- \mathrm{i} n
    t^\ast}.
  \label{eq:2.16}
\end{equation}
Since the coefficients, $\hat{\mu}[n]$ for $n = -N/2, \cdots, N/2-1$,
can be computed readily using the Fast Fourier Transform, we introduce
the truncated sum as an approximation in \eqref{eq:2.16}. The decay of
$\hat{\mu}[n]$ controls the error of this approximation. Therefore,
choosing $N$ to accurately solve the boundary integral equation yields
a spectrally accurate approximation of the double-layer potential.

For general problems, we do not know $\hat{K}[n]$ explicitly as we do
here. Instead, we compute an asymptotic expansion of
  the sharply peaked kernel. We determine the Fourier
coefficients of this asymptotic expansion explicitly.  Hence, we
evaluate its contribution to the double-layer potential using an
approximation like the one in \eqref{eq:2.16}. By removing the
kernel's sharp peak in this way, we are left with a smooth function to
integrate using the PTR$_{N}$. We present this
method to evaluate the double-layer potential for the
interior Dirichlet problem for Laplace's equation in Section
\ref{sec:doublelayer}, and the single-layer potential for the exterior
Neumann problem for Laplace's equation in Section
\ref{sec:singlelayer}.

\section{Double-layer potential for the interior Dirichlet problem}
\label{sec:doublelayer}

Consider a simply connected, open set denoted by
$D \subset \mathbb{R}^2$, with analytic boundary $\partial D$. Let
$\overline{D} = D \cup \partial D$. The function
$u \in C^2(D) \cap C^1(\overline{D})$ satisfies
\begin{subequations}
  \begin{gather}
    \Delta u = 0 \quad \text{in $D$}, \label{eq:3.1a}\\
    u = f \quad \text{on $\partial D$}, \label{eq:3.1b}
  \end{gather}
  \label{eq:3.1}
\end{subequations}
with $f$ an analytic function. We seek $u$ as the
double-layer potential~\cite{kress1999linear},
\begin{equation}
  u(\mathbf{x}) = \frac{1}{2\pi} \int_{\partial D}
  K(\mathbf{x},\mathbf{y})\mu(\mathbf{y})\mathrm{d}\sigma_{y}, \quad
  \mathbf{x} \in D,
  \label{eq:DLP}
\end{equation}
with
\begin{equation}
  K(\mathbf{x},\mathbf{y})=\mathbf{n}_{y} \cdot \frac{\mathbf{x} -
    \mathbf{y}}{| \mathbf{x} - \mathbf{y} |^{2}}.
  \label{eq:DLP-kernel}
\end{equation}
The density, $\mu$, satisfies the boundary integral equation,
\begin{equation}
  - \frac{1}{2} \mu(\mathbf{y}) + \frac{1}{2\pi} \int_{\partial D}
  K(\mathbf{y},\mathbf{y'})\mu(\mathbf{y}') \mathrm{d}\sigma_{y'} =
  f(\mathbf{y}), \quad \mathbf{y} \in \partial D.
  \label{eq:DLP-BIE}
\end{equation}
In what follows, we assume that we have solved \eqref{eq:DLP-BIE}
using PTR$_{N}$.

\begin{figure}[t]
  \centering
  \def\svgwidth{0.4\columnwidth} 
  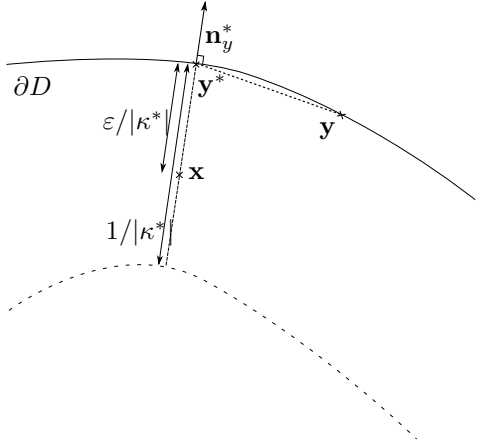
  \caption{Sketch of the quantities introduced in
    \eqref{eq:DLP-target} to study evaluation points close to the
    boundary.}
  \label{fig:sketch}
\end{figure}

To evaluate \eqref{eq:DLP} when $\mathbf{x}$ is close to the boundary,
we set
\begin{equation}
  \mathbf{x} = \mathbf{y}^{\ast} - \frac{\epsilon}
{  |\kappa^{\ast}| }\mathbf{n}_{y}^{\ast},
  \label{eq:DLP-target}
\end{equation}
where  $\mathbf{y}^{\ast}$ is the closest point to $\mathbf{x}$ on
the boundary, $\mathbf{n}_{y}^{\ast}$ is the unit, outward
normal at $\mathbf{y}^{\ast}$, $\kappa^{\ast}$ is the signed
curvature at $\mathbf{y}^{\ast}$, and $\epsilon >0$ is a small
parameter.  Fig.~\ref{fig:sketch} gives a sketch of these
quantities.  Substituting \eqref{eq:DLP-target} into
\eqref{eq:DLP-kernel} yields
\begin{equation}
  K \left(\mathbf{y}^{\ast} - \frac{\epsilon}
{  |\kappa^{\ast}| }\mathbf{n}_{y}^{\ast},
  \mathbf{y} \right) =
  \frac{|\kappa^{\ast}|}{\epsilon} \frac{\mathbf{n}_{y}
    \cdot |\kappa^{\ast}| ( \mathbf{y}^{\ast} - \mathbf{y} )/\epsilon -
    \mathbf{n}_{y} \cdot \mathbf{n}_{y}^{\ast}}{| \kappa^{\ast} (
    \mathbf{y}^{\ast} - \mathbf{y} )/\epsilon |^{2} - 2
    \mathbf{n}_{y}^{\ast} \cdot |\kappa^{\ast}| ( \mathbf{y}^{\ast} -
    \mathbf{y} )/\epsilon + 1}.
  \label{eq:asympt-kernel}
\end{equation}
We have written $K$ in \eqref{eq:asympt-kernel} to reveal its inherent
dependence on the stretched variable,
$\mathbf{y} = \mathbf{y}^{\ast} + \epsilon
\mathbf{Y}/|\kappa^{\ast}|$.

\subsection{Matched asymptotic expansion of the kernel}
\label{ssec:asymp}

We determine the matched asymptotic expansion of
\eqref{eq:asympt-kernel} \cite{bender1999advanced}.  Consider first
the outer expansion in which $\mathbf{y}^{\ast}$ and $\mathbf{y}$ are
held fixed and $\epsilon \to 0^{+}$, so that
$|\mathbf{Y}| \to \infty$. To leading order, we find that
\begin{equation}
  K^{\text{out}} \sim -\frac{|\kappa^{\ast}|}{\epsilon}
  \frac{\mathbf{n}_{y} \cdot \mathbf{Y}}{|\mathbf{Y}|^{2}} =
  \frac{\mathbf{n}_{y} \cdot ( \mathbf{y}^{\ast} - \mathbf{y} )}{|
    \mathbf{y}^{\ast} - \mathbf{y} |^{2} }.
  \label{eq:DLP-outer}
\end{equation}
The error of \eqref{eq:DLP-outer} is $O(\epsilon)$.  Since this outer
expansion 
is the kernel in \eqref{eq:DLP-BIE}, we find that
\begin{equation}
  \frac{1}{2\pi} \int_{\partial D} K^{\text{out}}(\mathbf{y}^{\ast} -
  \mathbf{y}) \mu(\mathbf{y}) \mathrm{d}\sigma_{y} =
  f(\mathbf{y}^{\ast}) + \frac{1}{2} \mu(\mathbf{y}^{\ast}).
  \label{eq:outerDLP}
\end{equation}

The inner expansion is \eqref{eq:asympt-kernel}
  written in terms of the stretched variable, $\mathbf{Y}$. We seek
an explicit parameterization of this inner expansion
using $ \mathbf{y}(t) = ( y_{1}(t), y_{2}(t) )$, with
$t \in [0,2\pi]$.  It follows that
$\mathrm{d}\sigma_{y} = | \mathbf{y}'(t) | \mathrm{d}t$, the unit
tangent is
$\boldsymbol{\tau}_{y}(t) = (y_{1}'(t),y_{2}'(t))/|\mathbf{y}'(t)|$,
the outward unit normal is
$\mathbf{n}_{y}(t) = (y_{2}'(t),-y_{1}'(t))/|\mathbf{y}'(t)|$, and the
signed curvature is
$ \kappa(t) = ( y_{1}'(t) y_{2}''(t) - y_{1}''(t)
y_{2}'(t))/|\mathbf{y}'({t})|^{3}$.
Let $\mathbf{y}^{\ast} = \mathbf{y}(t^{\ast})$ and
$\kappa^{\ast} = \kappa(t^{\ast})$ with $t^{\ast} \in [0,2\pi]$. We
introduce the stretched parameter, $t = t^{\ast} + \epsilon T$, and
find by expanding about $\epsilon = 0$ that
\begin{equation}
  \mathbf{y}(t^{\ast} + \epsilon T)  = \mathbf{y}(t^{\ast}) +
  \epsilon T | \mathbf{y}'(t^{\ast}) | \boldsymbol{\tau}_{y}(t^{\ast})
  - \frac{1}{2} \epsilon^{2} T^{2} \left[ \kappa^{\ast} |
    \mathbf{y}'(t^{\ast}) |^{2} \mathbf{n}_{y}(t^{\ast}) - (
    \boldsymbol{\tau}_{y}(t^{\ast}) \cdot \mathbf{y}''(t^{\ast}) )
    \boldsymbol{\tau}_{y}(t^{\ast}) \right] + O(\epsilon^{3}).
  \label{eq:3.10}
\end{equation}
It follows that
\begin{equation}
  \mathbf{n}_{y}(t^{\ast} + \epsilon T) \cdot |\kappa^{\ast}| (
  \mathbf{y}(t^{\ast}) -
  \mathbf{y}(t^{\ast} + \epsilon T ) ) = -\frac{1}{2} \epsilon^{2} T^{2}
  \gamma^{\ast} + O(\epsilon^{3}),
  \label{eq:3.11}
\end{equation}
\begin{equation}
  \mathbf{n}_{y}(t^{\ast}) \cdot |\kappa^{\ast}| ( \mathbf{y}(t^{\ast}) -
  \mathbf{y}(t^{\ast} + \epsilon T ) ) = \frac{1}{2} \epsilon^{2} T^{2}
  \gamma^{\ast} + O(\epsilon^{3}),
  \label{eq:3.12}
\end{equation}
\begin{equation}
  \mathbf{n}_{y}(t^{\ast} + \epsilon T ) \cdot
  \mathbf{n}_{y}(t^{\ast}) = 1 - \frac{1}{2} \epsilon^{2} T^{2}
  |\gamma^{\ast}|^{2} + O(\epsilon^{3}),
  \label{eq:3.13}
\end{equation}
and
\begin{equation}
  | \kappa(t^{\ast}) [ \mathbf{y}(t^{\ast}) - \mathbf{y}(t^{\ast} +
  \epsilon T ) ] |^{2} = \epsilon^{2} T^{2} | \gamma^{\ast} | +
  O(\epsilon^{3}),
  \label{eq:3.14}
\end{equation}
with
$\gamma^{\ast} = \text{sgn}[\kappa^{\ast}] | \kappa^{\ast}
\mathbf{y}'(t^{\ast}) |^{2}$. Here, $\text{sgn}[x] = x/|x|$ for
$x \ne 0$ and $\text{sgn}[x] = 0$ for $x = 0$. Substituting
\eqref{eq:3.11} -- \eqref{eq:3.14} into \eqref{eq:asympt-kernel}, we
find that
\begin{equation}
  K^{\text{in}}(T;\epsilon) = |\kappa(t^{\ast})| \frac{ - \epsilon
    -\frac{1}{2} \epsilon^{2} T^{2} \gamma^{\ast} +
    O(\epsilon^{3})}{\epsilon^{2} T^{2} | \gamma^{\ast} | +
    \epsilon^{2} + O(\epsilon^{3})}.
  \label{eq:3.15}
\end{equation}
Next, we substitute
$\epsilon^{2} T^{2} \sim 2 - 2 \cos( t - t^{\ast} )$ into
\eqref{eq:3.15} and determine that the leading order asymptotic
behavior of $K^{\text{in}}$ is given by
\begin{equation}
  K^{\text{in}}(t - t^{\ast};\epsilon) \sim |\kappa(t^{\ast})| \frac{ -
    ( \gamma^{\ast} + \epsilon ) + \gamma^{\ast} \cos(t -
    t^{\ast})}{(2 | \gamma^{\ast} | + \epsilon^{2} ) - 2 |
    \gamma^{\ast} | \cos(t - t^{\ast})}.
  \label{eq:DLP-inner}
\end{equation}
%
%

  The error of \eqref{eq:DLP-inner} is at most $O(\epsilon)$.  

\begin{figure}[t]
  \centering
   \includegraphics[width=0.40\linewidth]{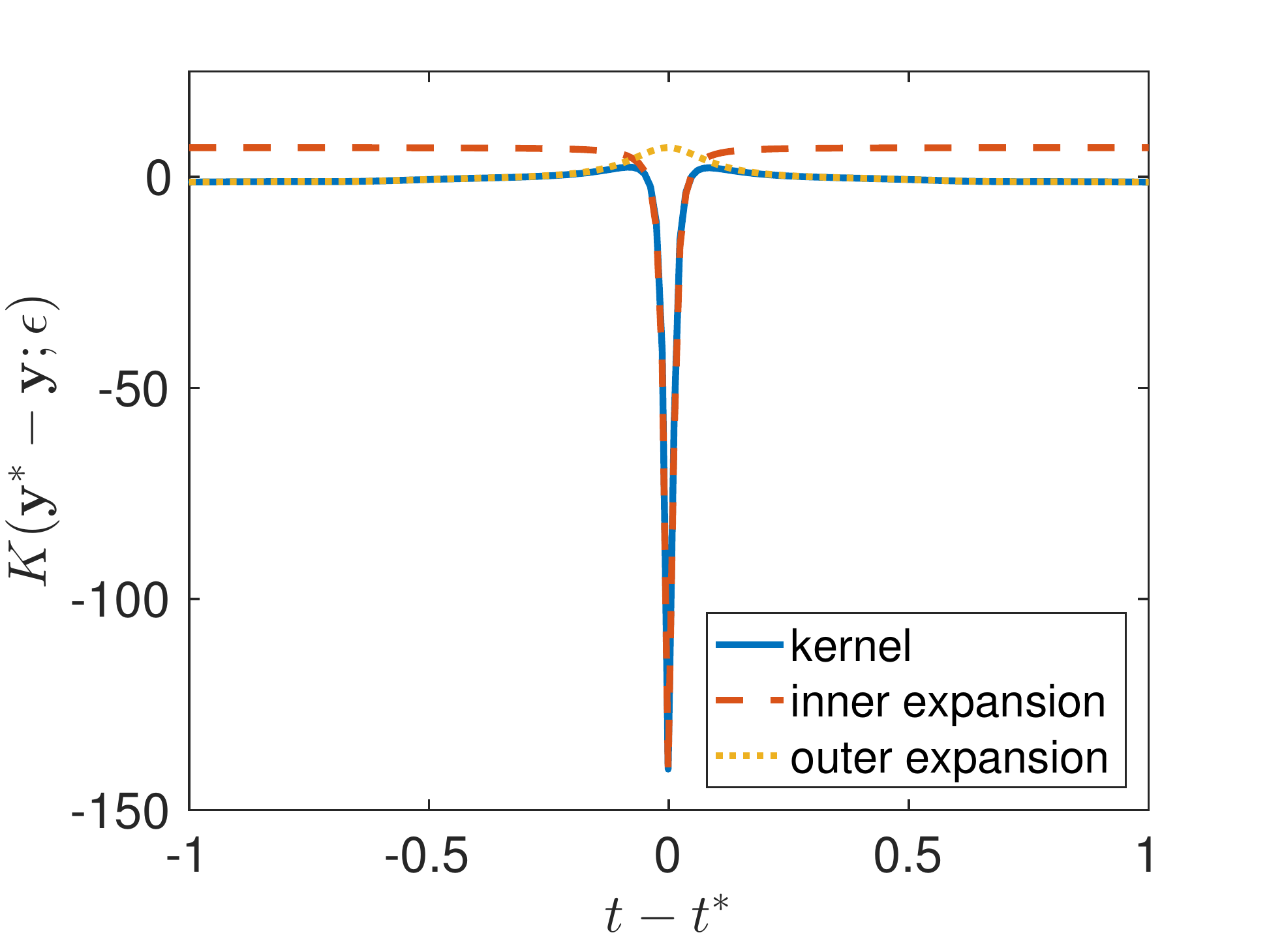}
  \includegraphics[width=0.40\linewidth]{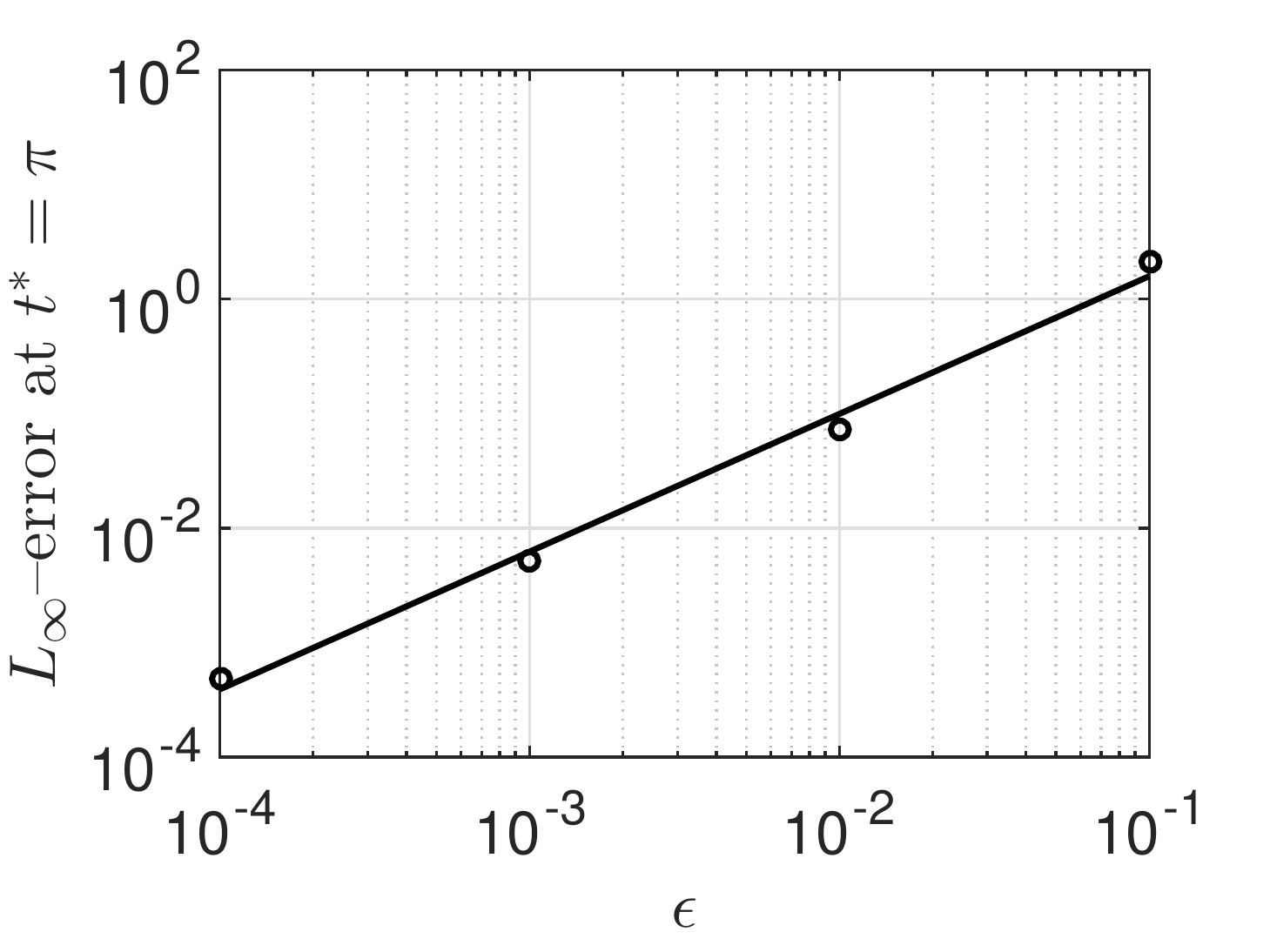}
  \caption{[Left] Plot of the kernel, $K$, given in
    \eqref{eq:asympt-kernel} (solid curve) and the leading order
    behavior of its inner expansion, $K^{\text{in}}$ given in
    \eqref{eq:DLP-inner} (dashed curve) and its outer expansion,
    $K^{\text{out}}$ given in \eqref{eq:DLP-outer} (dotted curve) as a
    function of $t - t^{\ast}$ with $t^{\ast} = \pi$ and
    $\epsilon = 0.1$ for the boundary curve,
    $r(t) = 1 + 0.3 \cos 5 t$. [Right] $L_{\infty}$--error of the
    matched asymptotic expansion given in \eqref{eq:DLP-matched}
    evaluated at $t^{\ast} = \pi$ for $\epsilon = 0.0001$, $0.001$,
    $0.01$, and $0.1$. These computed errors are plotted as
    circles. The solid curve gives the result of fitting this data to
    the function, $C \epsilon^{p}$.  This fit produced $p = 1.2034$
    indicating the $O(\epsilon)$ error of the matched asymptotic
    expansion.}
  \label{fig:4}
\end{figure}

To form the leading order matched asymptotic expansion,
we establish asymptotic matching in the overlap region
  of the outer and inner expansions. We first evaluate
\eqref{eq:DLP-outer} in the limit as
$\mathbf{y} \to \mathbf{y}^{\ast}$ and find that
\begin{equation}
  K^{\text{out}}  \to -
  \frac{\kappa^{\ast}}{2}, \quad \mathbf{y} \to \mathbf{y}^{\ast}.
\end{equation}
Next, we evaluate \eqref{eq:DLP-inner} in the limit as
$\epsilon \to 0^{+}$ and find that
\begin{equation}
  K^{\text{in}} \to - \frac{\kappa^{\ast}}{2}, \quad
  \epsilon \to 0^{+}.
\end{equation}
Thus, the overlapping value is $-\kappa^{\ast}/2$. It
  follows that the matched asymptotic expansion for the kernel of the
double-layer potential is given by
\begin{equation}
  K = K^{\text{out}} + K^{\text{in}} + \frac{\kappa^{\ast}}{2} + O(\epsilon),
  \quad \epsilon \to 0^{+},
  \label{eq:DLP-matched}
\end{equation}
with $K^{\text{out}}$ given in \eqref{eq:DLP-outer} and
$K^{\text{in}}$ given in \eqref{eq:DLP-inner}.  The error of this
matched asymptotic expansion has $O(\epsilon)$ error because the error
of $K^{\text{out}}$ given by \eqref{eq:DLP-outer} is $O(\epsilon)$,
and $K^{\text{in}}$ given by \eqref{eq:DLP-inner} is at most
$O(\epsilon)$.  For example, we plot $K$, $K^{\text{in}}$, and $K^{\text{out}}$ in
Fig.~\ref{fig:4} as a function of $t - t^{\ast}$ with $t^{\ast} = \pi$
and $\epsilon = 0.1$ for the boundary curve $r(t) = 1 + 0.3 \cos 5
t$. The right plot shows the $L_{\infty}$-error
made by \eqref{eq:DLP-matched} as a function of $ \epsilon$.
The solid curve is the linear fit through these data and has slope
$1.2034$ consistent with the $O(\epsilon)$ error.

\subsection{Fourier coefficients of $K^{\text{in}}$}

The inner expansion given by \eqref{eq:DLP-inner} accurately captures
the sharp peak of the kernel at $t = t^{\ast}$ in the
limit as $\epsilon \to 0^{+}$ as shown in
  Fig.~\ref{fig:4}. To avoid using PTR$_{N}$ to
integrate over this sharp peak, we seek the Fourier coefficients,
\begin{equation}
  \hat{K}^{\text{in}}[n] = \frac{1}{2\pi} \int_{0}^{2\pi}
  K^{\text{in}}(t;\epsilon) e^{-\mathrm{i} n t} \mathrm{d}t,
  \label{eq:3.19}
\end{equation}
so that we may use an approximation similar to that given in
\eqref{eq:2.16}. To do so, we rewrite \eqref{eq:DLP-inner} as
\begin{equation}
  K^{\text{in}}(t - t^{\ast};\epsilon) = - \frac{| \kappa(t^{\ast})
    |}{C_{0}} \frac{ \frac{1}{2} A_{0} - A_{1} \cos(t - t^{\ast})}{ 1 +
      C_{1} \cos(t - t^{\ast} )},
    \label{eq:DLP-FP}
\end{equation}
with $A_{0} = 2 ( \gamma^{\ast} + \epsilon ),$
$A_{1} = \gamma^{\ast},$ $C_{0} = 2 | \gamma^{\ast} | + \epsilon^{2},$
and $C_{1} = - 2 |\gamma^{\ast}| / C_{0}.$ Equation \eqref{eq:DLP-FP}
gives $K^{\text{in}}$ as a rational function of trigonometric
polynomials which have been studied by Geer~\cite{geer1995rational} in
the context of constructing Fourier-Pad\'e approximations. Since
$|C_{1}| < 1$, we have
 \begin{equation}
  \hat{K}^{\text{in}}[n] = 
  \begin{cases}
    \frac{1 + \rho^{2}}{1 - \rho^{2}} \left( \frac{A_{0}}{2} + A_{1}
      \rho \right), & n = 0\\
    \frac{1 + \rho^{2}}{1 - \rho^{2}} \left( \frac{A_{0}
        \rho^{|n|}}{4} + A_{1} ( \rho^{|n|-1} + \rho^{|n|+1} )
    \right), & n \neq 0
    \end{cases}, 
    \label{eq:DLP-Fourier}
\end{equation} 
where $\rho = \left(\sqrt{1 - C_{1}^{2}} - 1 \right)/C_{1}$.

We find that we can improve on our approximation by considering the
specific case in which the boundary is a circle of radius $a$. For
that case, $K^{\text{out}} = - \kappa^{\ast}/2$ which cancels with the
asymptotic matching term in \eqref{eq:DLP-matched}. If we set
\begin{subequations}
  \begin{align}
    A_{0} &= 2 ( \gamma^{\ast} + \epsilon - \epsilon | \gamma^{\ast} |),\\
    A_{1} &= \gamma^{\ast} - \epsilon | \gamma^{\ast} |,\\
    C_{0} &= 2 ( | \gamma^{\ast} | - \epsilon \gamma^{\ast} ) +
            \epsilon^{2},\\
    C_{1} &= - 2 ( | \gamma^{\ast} | - \epsilon \gamma^{\ast} ) / C_{0},
  \end{align}
  \label{eq:DLP-FPcoeffs}
\end{subequations}
instead of the coefficients defined above, we find that
\eqref{eq:DLP-FP} gives the exact evaluation of the kernel at
$r = a ( 1 - \epsilon )$. For this reason, we use
\eqref{eq:DLP-FPcoeffs} in \eqref{eq:DLP-FP} and
\eqref{eq:DLP-Fourier} in practice. These coefficients just include
the $O(\epsilon^3T^2)$ terms in the asymptotic
expansion of $K^{\text{in}}$ for a general boundary.

To compute the contribution by $K^{\text{in}}$ to the
  double-layer potential, we use the approximation
\begin{equation}
  \frac{1}{2\pi} \int_{0}^{2\pi} K^{\text{in}}(t - t^{\ast};\epsilon)
  \mu(t) | \mathbf{y}'(t) | \mathrm{d}t \approx \sum_{n = -N/2}^{N/2-1}
  \hat{K}^{\text{in} \ast}[n] \hat{\mu}_y[n] e^{-\mathrm{i}n t^\ast},
  \label{eq:innerDLP}
\end{equation}
with
\begin{equation}
  \hat{\mu}_y[n] = \frac{1}{2\pi} \int_{0}^{2\pi} \mu(t) | \mathbf{y}'(t)
  | e^{-\mathrm{i} n t} \mathrm{d}t.
  \label{eq:DLP-muhat}
\end{equation}
We use \eqref{eq:DLP-Fourier} and compute \eqref{eq:DLP-muhat} using
the Fast Fourier Transform to evaluate the approximation in
\eqref{eq:innerDLP}. Provided that $N$ is chosen to solve boundary
integral equation \eqref{eq:DLP-BIE} so that $\mu(t) | \mathbf{y}'(t)
|$ is sufficiently resolved, the approximation in \eqref{eq:innerDLP}
is spectrally accurate.

\subsection{Evaluating the double-layer potential}

The new method developed here for close evaluation of the double-layer
potential uses \eqref{eq:DLP-matched}.  For convenience, let us
introduce the residual kernel,
\begin{equation}
  \tilde{K} = K - K^{\text{out}} - K^{\text{in}} -
  \frac{\kappa^{\ast}}{2}.
  \label{eq:3.27}
\end{equation}
$\tilde{K} = O(\epsilon)$ and more importantly, it does not have a
sharp peak about $t = t^{\ast}$.  We rewrite the double-layer
potential as
\begin{multline}
  u\left(\mathbf{y}^{\ast} - \frac{\epsilon}
{  |\kappa^{\ast}| }\mathbf{n}_{y}^{\ast} \right) = \frac{1}{2\pi}
  \int_{0}^{2\pi} \tilde{K}(t - t^{\ast};\epsilon) \mu(t) |
  \mathbf{y}'(t) | \mathrm{d}t + \frac{1}{2\pi} \int_{0}^{2\pi}
  K^{\text{out}}(t - t^{\ast};\epsilon)
  \mu(t) | \mathbf{y}'(t) | \mathrm{d}t \\
  + \frac{1}{2\pi} \int_{0}^{2\pi} K^{\text{in}}(t -
  t^{\ast};\epsilon) \mu(t) | \mathbf{y}'(t) | \mathrm{d}t +
  \frac{\kappa^{\ast}}{4\pi} \int_{0}^{2\pi} \mu(t) |
  \mathbf{y}'(t) | \mathrm{d}t.
  \label{eq:3.26}
\end{multline}
Substituting \eqref{eq:outerDLP} and \eqref{eq:innerDLP} into
\eqref{eq:3.26}, we obtain
\begin{equation}
  u\left(\mathbf{y}^{\ast} - \frac{\epsilon}
{  |\kappa^{\ast}| }\mathbf{n}_{y}^{\ast} \right) \approx \frac{1}{2\pi}
  \int_{0}^{2\pi} \left[ \tilde{K}(t - t^{\ast};\epsilon) +
    \frac{\kappa^{\ast}}{2} \right] \mu(t) | \mathbf{y}'(t) |
  \mathrm{d}t + f(\mathbf{y}(t^{\ast})) + \frac{1}{2} \mu(t^{\ast})
  + \sum_{n = -N/2}^{N/2-1} \hat{K}^{\text{in} \ast}[n]
  \hat{\mu}_y[n] e^{- \mathrm{i} n t^\ast}.
  \label{eq:3.29}
\end{equation}
Applying PTR$_{N}$ with $t_{j} = 2\pi j/N$ to the
remaining integral in \eqref{eq:3.29}, we arrive at
\begin{equation}
 u\left(\mathbf{y}^{\ast} - \frac{\epsilon}
{  |\kappa^{\ast}| }\mathbf{n}_{y}^{\ast} \right)\approx \frac{1}{N} \sum_{j =
    1}^{N} \left[ \tilde{K}( t_{j} - t^{\ast};\epsilon ) +
    \frac{\kappa^{\ast}}{2} \right] \mu(t_{j}) |\mathbf{y}'(t_{j})|
  + f(\mathbf{y}(t^{\ast})) + \frac{1}{2} \mu(t^{\ast}) + \sum_{n =
    -N/2}^{N/2-1} \hat{K}^{\text{in} \ast}[n] \hat{\mu}_y[n] e^{-
    \mathrm{i} n t^\ast}.
  \label{eq:DLP-asymptotic}
\end{equation}

Equation \eqref{eq:DLP-asymptotic} gives our method for computing the
double-layer potential for close evaluation points. It avoids aliasing
incurred by the sharp peak of $K^{\text{in}}$ by using
\eqref{eq:innerDLP}. Integration of $K^{\text{out}}$ 
  is replaced by
  $f(\mathbf{y}^{\ast}) + \frac{1}{2} \mu(\mathbf{y}^{\ast})$, which
  comes from evaluating boundary integral equation \eqref{eq:DLP-BIE}
  at $\mathbf{y}^{\ast}$.
PTR$_{N}$ is now used only to integrate the term with
the kernel, $\tilde{K} + \kappa^{\ast}/2$.
  This term is important for taking into account additional, non-local
  contributions to the double-layer potential, which may be
  significant.

\subsection{Numerical results}\label{ssec:DLPresults}

We present results of this method for evaluating the double-layer
potential by computing the harmonic function,
$u(\mathbf{x}) = -\frac{1}{2\pi} \log | \mathbf{x} - \mathbf{x}_{0} |$
with $\mathbf{x}_{0} ={ (1.85,1.65)}$.  We compute the solution
interior to the boundary curve, $r(t) = 1.55 + 0.4 \cos 5 t$.  The
Dirichlet data in \eqref{eq:3.1b} is determined by evaluating the
harmonic function on the boundary.  Boundary integral equation
\eqref{eq:DLP-BIE} is solved using PTR$_{N}$ and we
use the resulting density, $\mu(t_{j})$ with $t_{j} = 2\pi j/N$ for
$j = 1, \cdots, N$ in the double-layer potential. We evaluate the
double-layer potential using two methods: (1)
  PTR$_{N}$ and (2) asymptotic PTR$_{N}$, the new method given in
\eqref{eq:DLP-asymptotic}.  We present the solution on a body-fitted
grid in which evaluation points are found by moving
along the normal into the domain from boundary grid
points.  The solution is evaluated on a grid of $200$ equispaced
points along each normal starting at the boundary until we reach a
distance $1/\kappa_{\max}$ from the boundary, where
$\kappa_{\max} = \max_{0 \le t^{\ast} < 2 \pi} |\kappa(t^{\ast})|$.
This grid captures the boundary layer, but does not coincide exactly
with it since the boundary layer depends on the local curvature.
For regions of high curvature, this body-fitted grid
extends beyond the boundary layer.
  
\begin{figure}[ht!]
  \centering
  \includegraphics[width=0.75\linewidth]{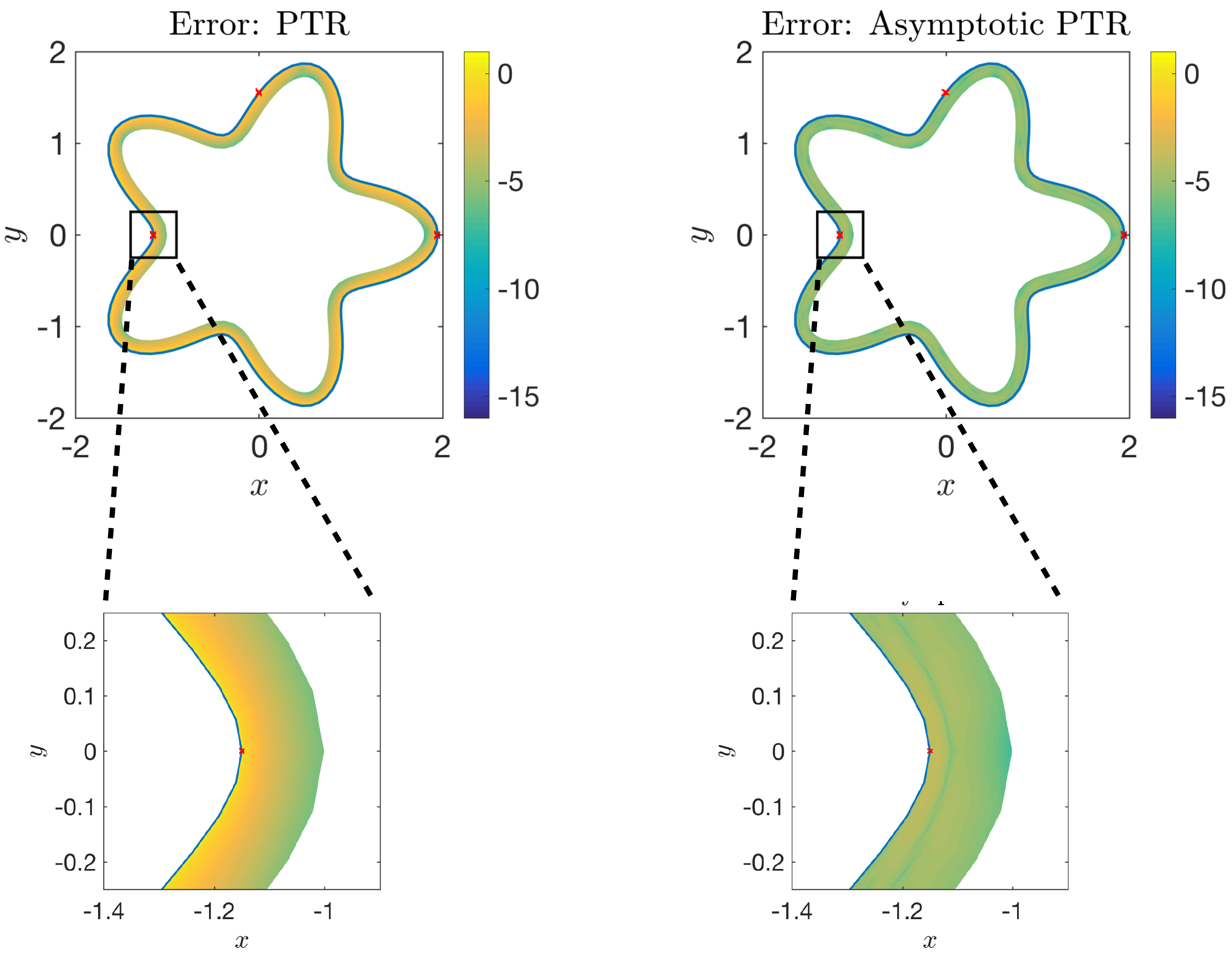}
  \caption{[Left] Plot of absolute error
      ($\log_{10}$-scale) in computing the double-layer potential
      using PTR$_{128}$ for the boundary
    $r(t) = {1.55 + 0.4 \cos 5 t }$ for the Dirichlet data,
    $f(\mathbf{y}) = \frac{1}{2\pi} \log | \mathbf{y} - \mathbf{x}_{0}
    |$
    with $\mathbf{x}_{0} = {(1.85, 1.65)}$.  [Right] Plot of
    absolute error ($\log_{10}$-scale) in computing
      the double-layer potential using the asymptotic PTR$_{128}$
      given in \eqref{eq:DLP-asymptotic} for the same problem. The
    ``$\times$'' symbols on the boundary indicates the points
    corresponding to $t^{\ast}=0$, $t^{\ast}=\pi/2$, and
    $t^{\ast}=\pi$. }
  \label{fig:6}
\end{figure}

\begin{figure}[ht!]
  \centering
  \includegraphics[width=0.31\linewidth]{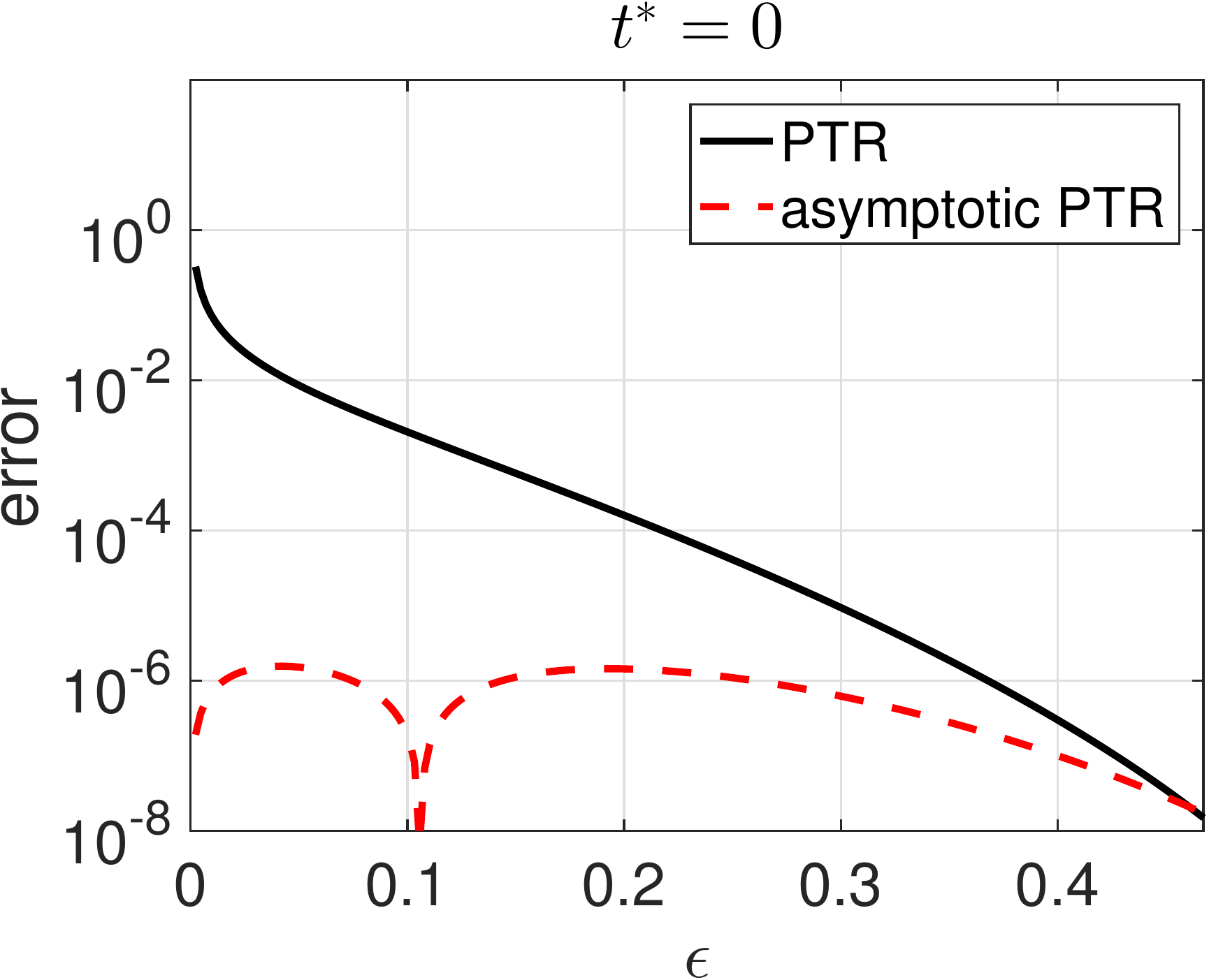} \, \,
  \includegraphics[width=0.31\linewidth]{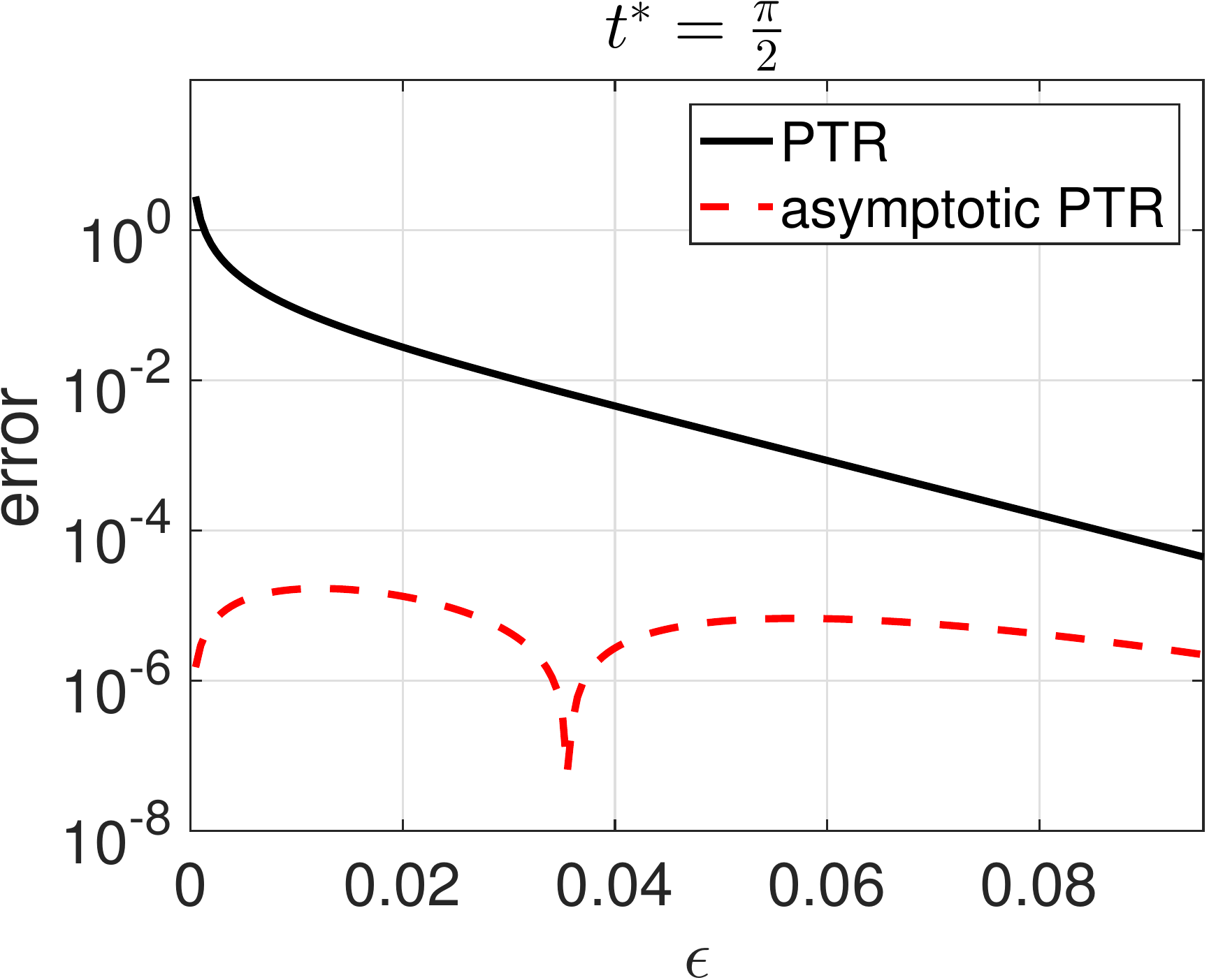} \, \,
  \includegraphics[width=0.31\linewidth]{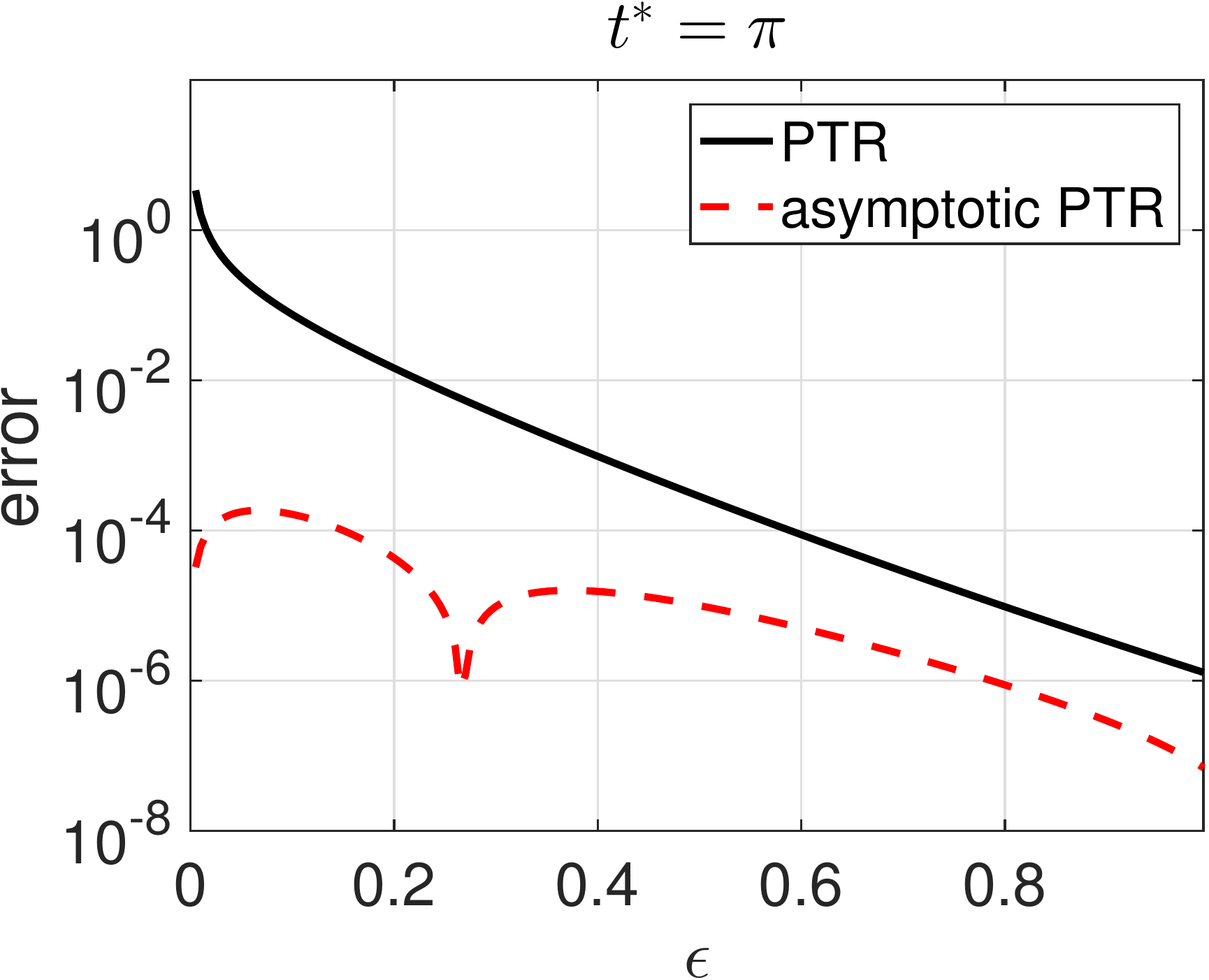}  
  \caption{Plot of the absolute error as a function of $\epsilon$ made
    in evaluating the double-layer potential for different $t^\ast$
    values using PTR$_{128}$ (solid curve) and using
    asymptotic PTR$_{128}$ (dashed curve).}
  \label{fig:7}
\end{figure}

In {Fig.}~\ref{fig:6} we show the errors
  ($\log_{10}$-scale) in computing the double-layer potential using
  PTR$_{128}$ and asymptotic PTR$_{128}$. These results show that
asymptotic PTR$_{N}$ produces errors that are several
orders of magnitude smaller than those of the
PTR$_{N}$. To give an indication of this
improvement, the $L_{\infty}$ error is $8.03$ for
PTR$_{128}$ and $1.85 \times 10^{-4}$ for asymptotic
PTR$_{128}$. To examine this error in more detail,
we show in {Fig.}~\ref{fig:7} a plot of the error in computing the
double-layer potential evaluated at the points indicated in
Fig.~\ref{fig:6} ($t^{\ast} = 0$, $t^{\ast} = \pi/2$, and
$t^{\ast} = \pi$) as a function of $\epsilon$.  These three cases are
plotted over different ranges of $\epsilon$ corresponding to
$0 < \epsilon < \kappa(t^{\ast})/\kappa_{\max}$.  We observe that
asymptotic PTR$_{N}$ does significantly better than
PTR$_{N}$ for small $\epsilon$ as expected. It
reduces the $O(1)$ error by at least four orders of magnitude.  We
find similar results over all values of
$t^{\ast}$.

We can further improve the new method using the identity for the
  double-layer potential~\cite{kress1999linear} (see
  \eqref{kernel_identities}) to rewrite \eqref{eq:DLP} as follows:
\begin{equation}
  u(\mathbf{x}) = \frac{1}{2\pi} \int_{\partial D}
  K(\mathbf{x},\mathbf{y})(\mu(\mathbf{y}) - \mu(\mathbf{y}^\ast)
  )\mathrm{d}\sigma_{y} - \mu(\mathbf{y}^\ast), \quad
  \mathbf{x} \in \partial D.
  \label{eq:DLP+subtraction}
\end{equation} 
In \eqref{eq:DLP+subtraction}, the integrand is now smoother as it
vanishes at the point $\mathbf{y} = \mathbf{y}^\ast$, and the error
using PTR$_{N}$ drastically decreases.  Applying
asymptotic PTR$_{N}$ to \eqref{eq:DLP+subtraction},
we obtain
\begin{equation}
  u\left(\mathbf{y}^{\ast} - \frac{\epsilon}
    {  |\kappa^{\ast}| }\mathbf{n}_{y}^{\ast} \right) \approx
    \frac{1}{N} \sum_{j = 1}^{N} \left[ \tilde{K}( t_{j} -
      t^{\ast};\epsilon ) + \frac{\kappa^{\ast}}{2} \right]
    (\mu(t_{j})-\mu(t^\ast)) |\mathbf{y}'(t_{j})| +
    f(\mathbf{y}(t^{\ast}))  + \sum_{n = -N/2}^{N/2-1}
    \hat{K}^{\text{in} \ast}[n] \hat{\mu}^\ast_y[n] e^{- \mathrm{i} n
      t^\ast}
  \label{eq:asymptoticDLP+subtraction}
\end{equation} 
with
\begin{equation}
  \hat{\mu}^\ast_y[n] = \frac{1}{2\pi} \int_{0}^{2\pi} (\mu(t)
  -\mu(t^\ast)) | \mathbf{y}'(t) | e^{-\mathrm{i} n t} \mathrm{d}t.
  \label{eq:DLP-muhatast}
\end{equation}
For the example problem discussed above, the $L_{\infty}$ error is
$7.13 \times 10^{-5}$ for PTR$_{128}$ applied to
\eqref{eq:DLP+subtraction}, and $2.54 \times 10^{-5}$ for asymptotic
PTR$_{128}$ given by
\eqref{eq:asymptoticDLP+subtraction}.
This additional improvement \eqref{eq:DLP+subtraction} is only valid for the 
double-layer potential and does not generalize.

\section{Single-layer potential for the exterior Neumann problem}
\label{sec:singlelayer}

We now consider the exterior Neumann problem,
\begin{subequations}
  \begin{gather}
    \Delta v = 0 \quad \text{in $\mathbb{R}^{2} \backslash
      \bar{D}$}, \label{eq:4.1a}\\
    \frac{\partial v}{\partial n} = g \quad \text{on $\partial
      D$}, \label{eq:4.1b}
  \end{gather}
  \label{eq:4.1}
\end{subequations}
with $g$ denoting an analytic function satisfying
\begin{equation}
  \int_{\partial D} g(\mathbf{y}) \mathrm{d}\sigma_{y} = 0.
  \label{eq:4.2}
\end{equation}
We seek $v$ as the single-layer potential~\cite{kress1999linear},
\begin{equation}
  v(\mathbf{x}) = \frac{1}{2\pi} \int_{\partial D}
  S(\mathbf{x},\mathbf{y}) \varphi(\mathbf{y}) \mathrm{d}\sigma_{y},
  \quad \mathbf{x} \in \mathbb{R}^{2} \backslash \bar{D},
  \label{eq:SLP}
\end{equation}
with
\begin{equation}
  S(\mathbf{x},\mathbf{y}) = - \log | \mathbf{x} - \mathbf{y} | .
  \label{eq:SLP-kernel}
\end{equation}
The density, $\varphi(\mathbf{y})$, satisfies the boundary integral
equation,
\begin{equation}
  - \frac{1}{2} \varphi(\mathbf{y}) + \frac{1}{2\pi} \int_{\partial D}
  \frac{\partial S(\mathbf{y},\mathbf{y'}) }{\partial n_{y} }
  \varphi(\mathbf{y}') \mathrm{d}\sigma_{y'} = g(\mathbf{y}), \quad
  \mathbf{y} \in \partial D.
  \label{eq:SLP-BIE}
\end{equation}

To study the close evaluation of \eqref{eq:SLP}, we now set 
\begin{equation}
  \mathbf{x} = \mathbf{y}^{\ast} + \frac{\epsilon}{|\kappa^{\ast}|}
  \mathbf{n}_{y}^{\ast}.
  \label{eq:SLP-target}
\end{equation}
Substituting \eqref{eq:SLP-target} into \eqref{eq:SLP-kernel}, we obtain
\begin{equation}
  S\left(\mathbf{y}^{\ast} + \frac{\epsilon}
{  |\kappa^{\ast}| }\mathbf{n}_{y}^{\ast},
  \mathbf{y} \right) = - \log \epsilon +
  \log | \kappa^{\ast} | - \frac{1}{2} \log\left(  | \kappa^{\ast}
    ( \mathbf{y}^{\ast} - \mathbf{y} )/\epsilon |^{2} + 2 
    \mathbf{n}_{y}^{\ast} \cdot | \kappa^{\ast} | ( \mathbf{y}^{\ast}
    - \mathbf{y} )/\epsilon + 1 \right).
  \label{eq:asympt-SLPkernel}
\end{equation}
Just as we have done for $K$ in \eqref{eq:asympt-kernel}, we have
written \eqref{eq:asympt-SLPkernel} to show the underlying dependence
on the stretched variable,
$\mathbf{y} = \mathbf{y}^{\ast} + \epsilon
\mathbf{Y}/|\kappa^{\ast}|$.

The outer expansion of \eqref{eq:asympt-SLPkernel} is
$S^{\text{out}} \sim - \log |\mathbf{y^\ast} -\mathbf{y}|$.  
  This outer expansion is singular. In contrast to the double-layer
  potential, this outer expansion does not correspond to the kernel
  boundary integral equation \eqref{eq:SLP-BIE}.  One could use a high
  order quadrature rule that explicitly takes into account this
  singularity~\cite{sidi1988quadrature, kress1991boundary}. However,
  we choose to not use one here because we find in the numerical
  examples below that our method significantly reduces the dominant error. 
To compute the inner expansion, $S^{\text{in}}$, we introduce the
stretched parameter, $t = t^{\ast} + \epsilon T$ into the same
parameterization of the boundary used for the double-layer
potential. Making use of \eqref{eq:3.11} and \eqref{eq:3.13}, we find
by expanding as $\epsilon \to 0^{+}$ that
\begin{equation}
  S^{\text{in}}(T;\epsilon) = \log | \kappa^{\ast} | - \frac{1}{2}
  \log\left( \epsilon^{2} T^{2} | \gamma^{\ast} | + \epsilon^{2} +
    O(\epsilon^{3}) \right).
  \label{eq:4.7}
\end{equation}
Substituting $\epsilon^{2} T^{2} \sim 2 - 2 \cos( t - t^{\ast} )$, we
find that to leading order,
\begin{equation}
  S^{\text{in}}(t - t^{\ast};\epsilon) \sim \log | \kappa^{\ast} | -
  \frac{1}{2} \log\left( ( 2 | \gamma^{\ast} | + \epsilon^{2} ) - 2 |
    \gamma^{\ast} | \cos( t - t^{\ast}) \right).
  \label{eq:4.8}
\end{equation}

\subsection{Fourier coefficients of $S^{\text{in}}$}
Using the modified coefficients introduced in \eqref{eq:DLP-FPcoeffs},
we write \eqref{eq:4.8} as
\begin{equation}
  S^{\text{in}}(t - t^{\ast};\epsilon) \sim \log |\kappa^{\ast}| -
  \frac{1}{2} \log C_{0}  - \frac{1}{2} \log[ 1 + C_{1} \cos(t -
  t^{\ast}) ].
  \label{eq:4.10}
\end{equation}
We now seek to compute 
\begin{equation}
  \hat{S}^{\text{in}}[n] =  \delta_{n,0} \left[ \log |\kappa^{\ast}|
    - \frac{1}{2} \log C_{0} \right] - \frac{1}{2\pi} \int_{0}^{2\pi}
  \frac{1}{2} \log[ 1 + C_{1} \cos(t - t^{\ast}) ] e^{-\mathrm{i} n t}
  \mathrm{d}t,
  \label{eq:4.11}
\end{equation}
with $\delta_{n,0}$ denoting the Kronecker delta.  To compute the
integral in \eqref{eq:4.11}, we start with
\begin{equation}
  \frac{\mathrm{d}}{\mathrm{d}t} \log[ 1 +
  C_{1} \cos(t - t^{\ast}) ] = \frac{C_{1} \sin(t - t^{\ast})}{1 +
    C_{1} \cos(t - t^{\ast})}.
  \label{eq:4.12}
\end{equation}
The right-hand side of \eqref{eq:4.12} is another example of a
rational trigonometric function studied by
Geer~\cite{geer1995rational}. It can be readily shown that
\begin{equation}
  \frac{1}{2\pi} \int_{0}^{2\pi} \frac{\sin(t -
    t^{\ast})}{1 + C_{1} \cos(t - t^{\ast})} e^{\mathrm{i} n t}
  \mathrm{d}t = \text{sgn}(n) \frac{\mathrm{i}}{2} \frac{1 +
    \rho^{2}}{1 - \rho^{2}} \left( \rho^{|n| - 1} - \rho^{|n| + 1}
  \right),
\end{equation}
where $\rho = \left( \sqrt{1 - C_{1}^{2}} \right)/C_{1}$.  It follows
from term-by-term integration of the Fourier series with these
coefficients that
\begin{equation}
  \hat{S}^{\text{in}}[n] = \begin{cases}
    \displaystyle \log | \kappa^{\ast} | - \frac{1}{2} \log C_{0} -
    \frac{C_{1}}{2} \frac{1 + \rho^{2}}{1 - \rho^{2}} \left(
      \frac{1}{\rho} - \rho \right) \log( 1 - \rho) & n = 0,\\
    \displaystyle \frac{C_{1}}{4 |n|} \frac{1 +
      \rho^{2}}{1 - \rho^{2}} \left( \rho^{|n| - 1} - \rho^{|n| + 1}
    \right) & n \neq 0.
  \end{cases}
  \label{eq:4.14}
\end{equation}

\subsection{Evaluating the single-layer potential}

Given the inner expansion computed above, our method for evaluating
the single-layer potential is to compute an approximation of
\begin{equation}
  v\left(\mathbf{y}^{\ast} + \frac{\epsilon}
{  |\kappa^{\ast}| }\mathbf{n}_{y}^{\ast} \right) = \frac{1}{2\pi}
  \int_{0}^{2\pi} \tilde{S}(t - t^{\ast}) \varphi(t) | \mathbf{y}'(t) |
  \mathrm{d}t + \frac{1}{2\pi} \int_{0}^{2\pi} S^{\text{in}}(t -
  t^{\ast};\epsilon) \varphi(t) | \mathrm{y}'(t) | \mathrm{d}t,
  \label{eq:4.9}
\end{equation}
with $\tilde{S} = S - S^{\text{in}}$. We use
PTR$_{N}$ to evaluate the first integral with kernel,
$\tilde{S}$, and a truncated convolution sum to evaluate the second
integral with kernel, $S^{\text{in}}$,
\begin{equation}
  v\left(\mathbf{y}^{\ast} + \frac{\epsilon}
    {  |\kappa^{\ast}| }\mathbf{n}_{y}^{\ast} \right)  \approx
  \frac{1}{N} \sum_{j = 1}^{N} \tilde{S}(t_{j} - t^{\ast})
  \varphi(t_{j}) | \mathbf{y}'(t_{j}) | + \sum_{n = -N/2}^{N/2-1}
  \hat{S}^{\text{in}}[n] \hat{\varphi}_{y}[n] e^{-\mathrm{i} n t^\ast},
  \label{eq:4.15}
\end{equation}
where we compute
\begin{equation}
  \hat{\varphi}_{y}[n] = \frac{1}{2\pi} \int_{0}^{2\pi} \varphi(t) |
  \mathbf{y}'(t) | e^{-\mathrm{i} n t} \mathrm{d}t,
  \label{eq:4.16}
\end{equation}
using the Fast Fourier Transform.  Just as with the double-layer
potential, provided that $N$ is chosen so that it solves boundary
integral equation \eqref{eq:SLP-BIE} with sufficient accuracy, the
truncated convolution sum in \eqref{eq:4.15} is spectrally accurate.

\subsection{Numerical examples}\label{ssec:SLPresults}

We present results for the evaluation of the single-layer potential by
computing the harmonic function,
$v(\mathbf{x}) = ( \mathbf{x} - \mathbf{x}_{0} ) / | \mathbf{x} -
\mathbf{x}_{0} |^{2}$,
with $\mathbf{x}_0 = (0.1,0.4)$.  We compute the solution exterior to
the the boundary curve, $r(t) = 1.55 + 0.4 \cos 5 t$.  The Neumann
data in \eqref{eq:4.1b} is determined by computing the normal
derivative of the harmonic function on the boundary.  Boundary
integral equation \eqref{eq:SLP-BIE} is solved using
PTR$_{N}$ and we use the resulting density
$\varphi(t_{j})$ with $t_{j} = 2\pi j/N$ for $j = 1, \cdots, N$ in the
single-layer potential. We evaluate the single-layer potential using
two methods: (1) PTR$_{N}$ and (2) asymptotic
PTR$_{N}$, the new method given in \eqref{eq:4.15}.
We modify the body-fitted grid described above for the evaluation of
the double-layer potential to evaluate exterior points.  The solution
is evaluated on a grid of $200$ equispaced points along each normal
starting at the boundary until we reach a distance $1/\kappa_{\max}$
from the boundary.

\begin{figure}[h!]
  \centering
  \includegraphics[width=0.75\linewidth]{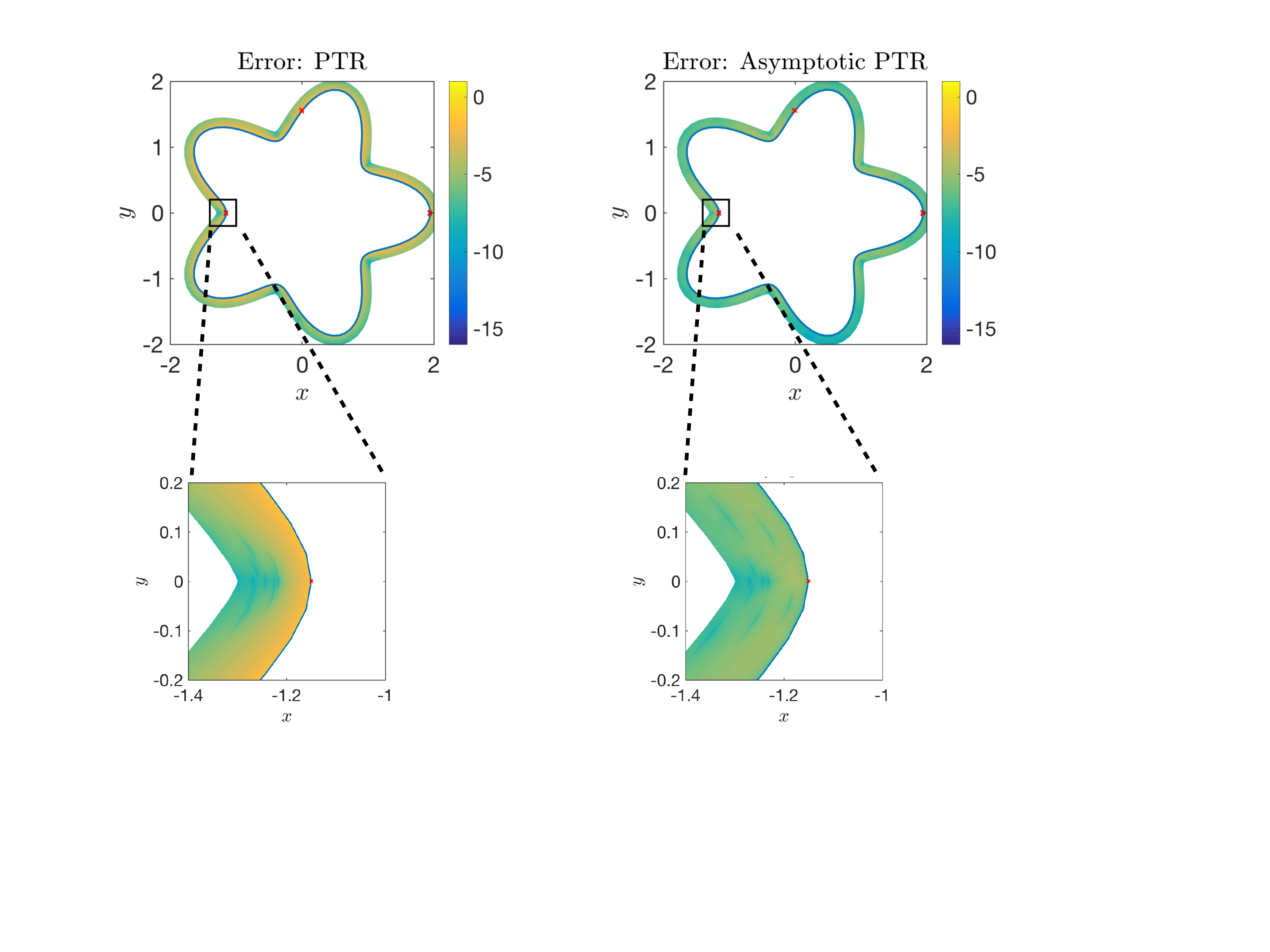}  
  \caption{[Left] Plot of the absolute error
      ($\log_{10}$-scale) in computing the single-layer potential
      using PTR$_{128}$ for the boundary
    $r(t) = {1.55 + 0.4 \cos 5 t }$ for the Neumann data,
    $g(\mathbf{y}) = \frac{\partial v}{\partial \mathbf{n}}$ with
    $v(\mathbf{x}) = ( \mathbf{x} - \mathbf{x}_{0} ) / | \mathbf{x} -
    \mathbf{x}_{0} |^{2}$,
    $\mathbf{x}_{0} = {(0.1, 0.4)}$.  [Right] Plot of the absolute
    error ($\log_{10}$-scale) in computing the
      single-layer potential using asymptotic PTR$_{128}$ given in
    \eqref{eq:4.15} for the same problem. The ``$\times$'' symbols on
    the boundary indicates the points corresponding to $t^{\ast}=0$,
    $t^{\ast}=\pi/2$ and $t^{\ast}=\pi$. }
  \label{fig:8}
\end{figure}

\begin{figure}[h!]
  \centering
  \includegraphics[width=0.31\linewidth]{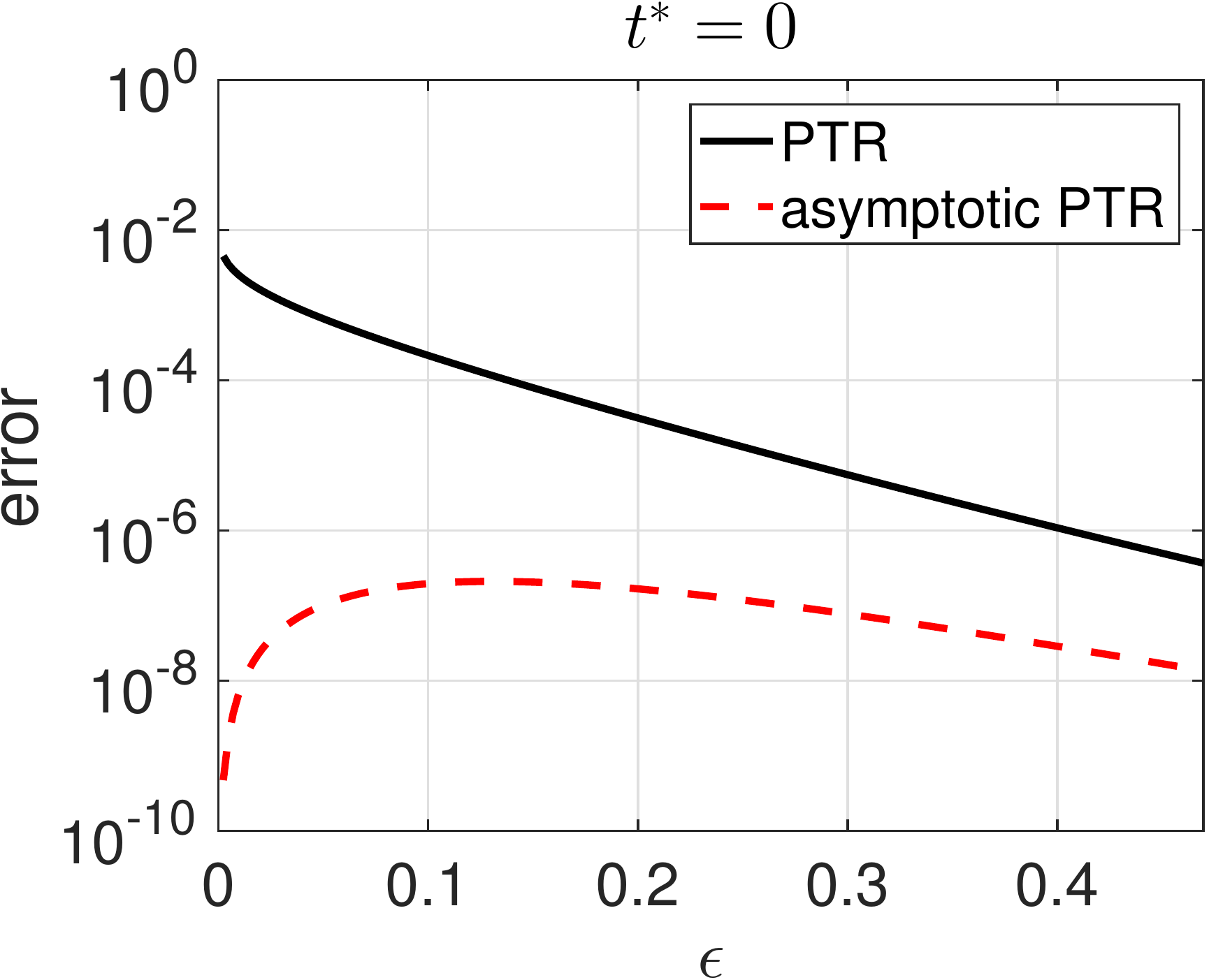} \, \, 
  \includegraphics[width=0.31\linewidth]{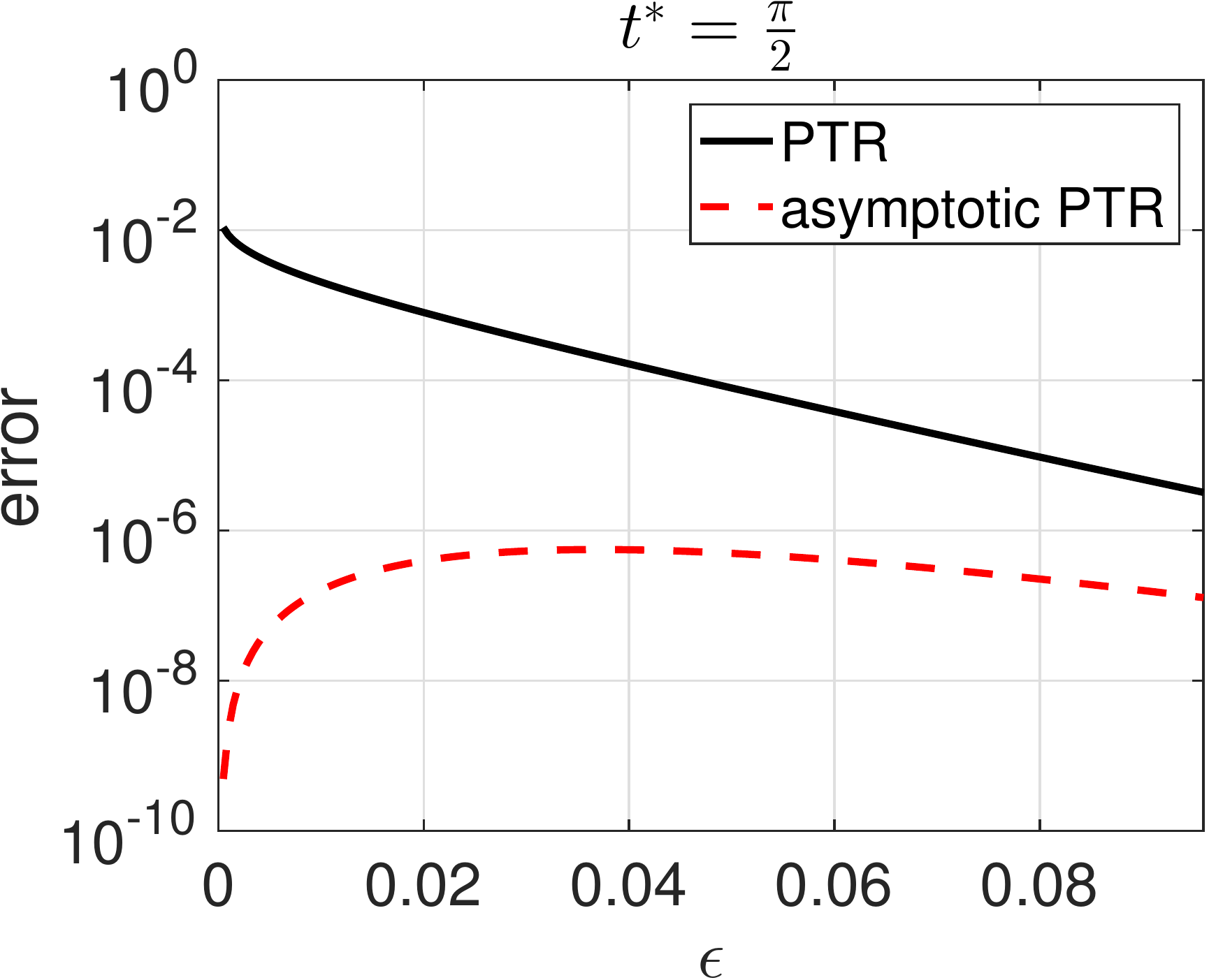} \, \, 
  \includegraphics[width=0.31\linewidth]{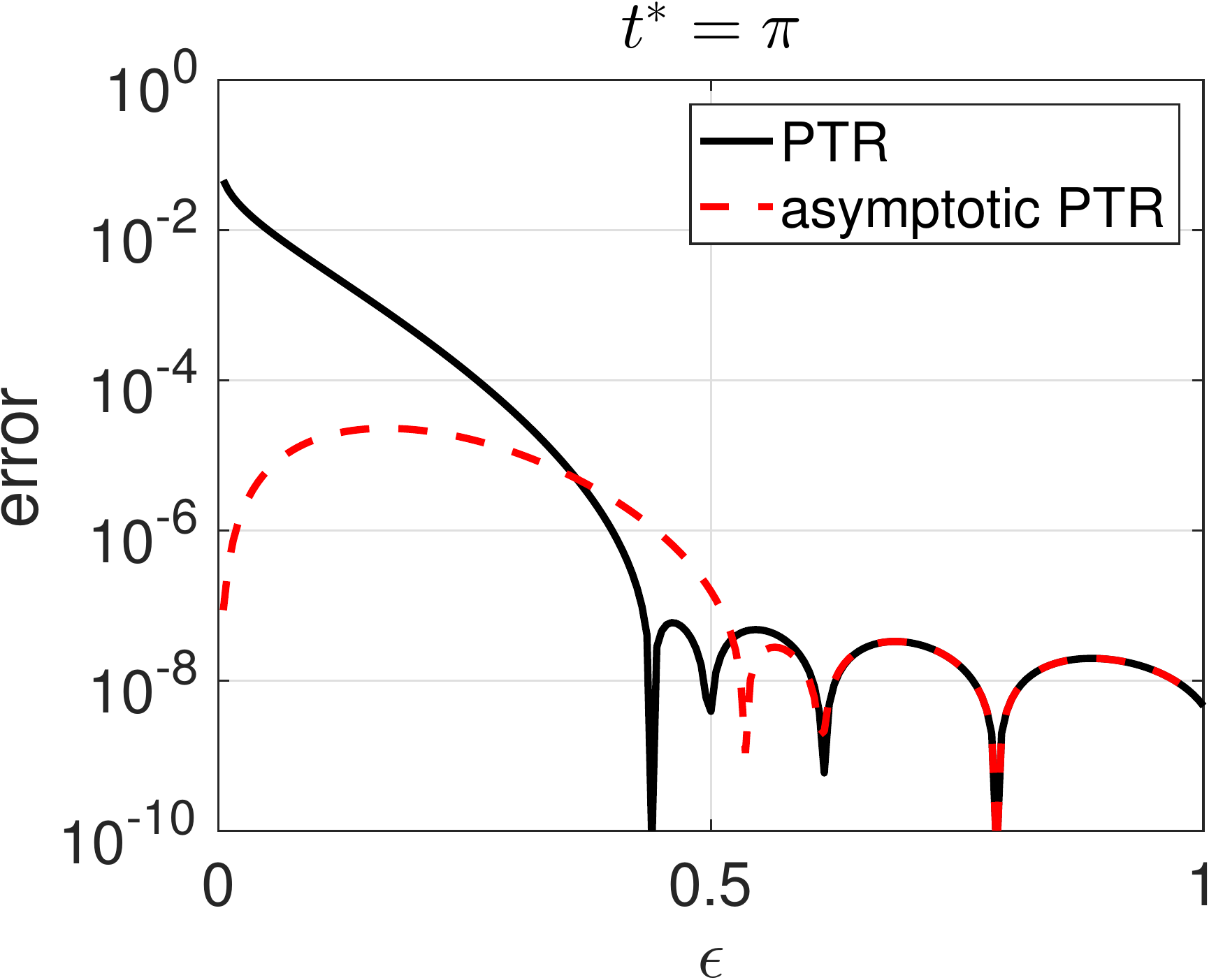} \, \, 
  \caption{Plot of the absolute error as a function of $\epsilon$ made
    in evaluating the single-layer potential for different $t^\ast$
    using PTR$_{128}$ (solid curve) and using
    asymptotic PTR$_{128}$ (dashed curve).}
    \label{fig:9}
\end{figure}

In {Fig.}~\ref{fig:8} we show the absolute error
  ($\log_{10}$-scale) in computing the single-layer potential using
the PTR$_{128}$ and using asymptotic
PTR$_{128}$.  The single-layer potential kernel is
not as sharply peaked as the double-layer potential kernel, so the
error in evaluating the single-layer potential is less than the error
when evaluating the double-layer potential. Even so, we still observe
a boundary layer of thickness $O(1/N)$ in which the error is $O(1)$
due to aliasing when using PTR$_{N}$. Asymptotic
PTR$_{N}$ effectively reduces the error in the
boundary layer. To give an indication of this improvement, the
$L_{\infty}$ error is $0.113$ for PTR$_{128}$ and
$5.39 \times 10^{-5}$ for asymptotic PTR$_{128}$.  In
Fig.~\ref{fig:9}, we plot the error in computing the single-layer
potential evaluated at the points indicated in Fig.~\ref{fig:8}
($t^{\ast} = 0$, $t^{\ast} = \pi/2$ and $t^{\ast} = \pi$) as a
function of $\epsilon$ with
$0 < \epsilon < \kappa(t^{\ast})/\kappa_{\max}$.  These plots show
that the asymptotic method reduces the error by at least 3 orders of
magnitude for small $\epsilon$. We find similar results over all
values of $t^{\ast}$.  For the case in which $t^{\ast} = \pi$, the
error of asymptotic PTR$_{N}$ becomes larger than
that for PTR$_{N}$ for $0.36 < \epsilon < 0.52$. For
this particular boundary curve, $\kappa_{\max}$ is attained at
$t^{\ast} = \pi$.  Hence, the body-fitted grid at $t^{\ast} = \pi$
plots the single-layer potential over $0 < \epsilon < 1$. For this
range of $\epsilon$, we consider points outside the boundary layer
where PTR$_{N}$ is competitive with, and may even
become more accurate than asymptotic PTR$_{N}$. In
fact, that is what is observed in Fig.~\ref{fig:9} for
$t^{\ast} = \pi$.

\section{General implementation}
\label{sec:implementation}

In the results presented above, we are using a body-fitted grid, in
which evaluation points are found by moving along the
  normal into the domain from boundary integration points.  Often the
solution is needed at more generally defined points.  Asymptotic
PTR$_{N}$ tacitly assumes in \eqref{eq:DLP-target}
for interior problems, or \eqref{eq:SLP-target} for exterior problems,
that $\mathbf{y}^{\ast}$, is the unique minimum distance from the
boundary to the evaluation point, $\mathbf{x}$.  In the examples, the
closest point on the boundary, $\mathbf{y}^{\ast}$ coincides with a
PTR$_{N}$ grid point from which we extended in the
normal direction.  We discuss here the more general problem.

Suppose we have an evaluation point in the domain, $\mathbf{x}$.  Then
the first problem to address is whether $\mathbf{x}$ is, in fact,
close enough to the boundary to warrant special attention and requires
use of asymptotic PTR$_{N}$. To solve this problem,
we make use of the identity for the double-layer
potential~\cite{kress1999linear},
\begin{equation}
  \frac{1}{2\pi}\int_{\partial D} \frac{\mathbf{n}_y \cdot (\mathbf{x}
    - \mathbf{y})}{|\mathbf{x} - \mathbf{y}|^2} \, \mathrm{d} \sigma_y
  = \begin{cases} 
    -1 & \mathbf{x} \in D\\
    -\frac{1}{2} & \mathbf{x} \in \partial D\\
    \,\,\, 0 & \mathbf{x} \in \mathbb{R}^2 \setminus \overline{D} 
  \end{cases}.
  \label{kernel_identities}
\end{equation}
Evaluating \eqref{kernel_identities} using PTR$_{N}$
suffers from the same aliasing problem that the more general layer
potential evaluations do. Thus, we use it to determine if $\mathbf{x}$
lies within the boundary layer. To do this, we set a user-defined
threshold for the error. If the error in evaluating
\eqref{kernel_identities} using PTR$_{N}$ is less
than this threshold value, we keep the result computed using
PTR$_{N}$.  Otherwise, we use the appropriate
asymptotic approximation.

To use these asymptotic approximations, we must determine the
parameter, $t^{\ast}$, where
$t^{\ast} = \min_{0 \le t < 2 \pi} | \mathbf{x} - \mathbf{y}(t)
|^{2}$.
For a general boundary, this problem may not have a unique
solution. In practice, we find a unique minimizer for evaluation
points that are identified to be in the boundary layer using
PTR$_{N}$ evaluation of \eqref{kernel_identities}.
Once $t^{\ast}$ is determined, all other quantities required for the
asymptotic approximations follow. Finally, we evaluate the solution of
the boundary value problem at hand at any point $\mathbf{x}$ using
either  PTR$_{N}$ or asymptotic
PTR$_{N}$.

\begin{figure}[h!]
  \centering
  \includegraphics[width=1\linewidth]{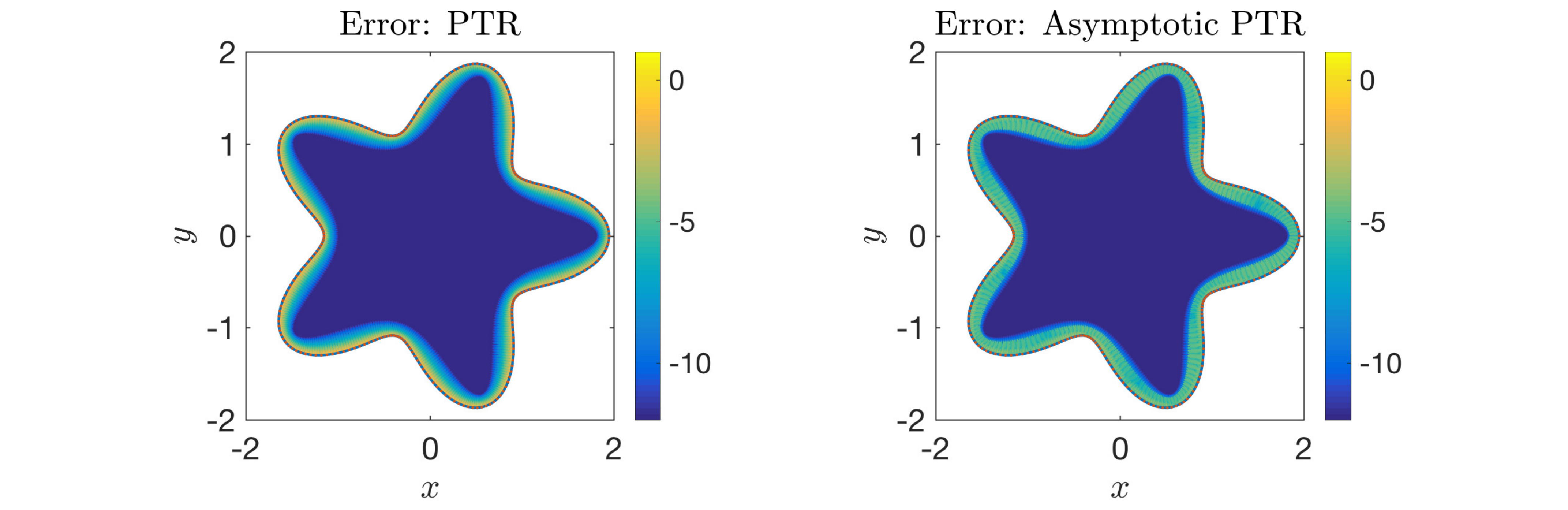}
  \caption{[Left] Plot of the absolute error
    ($\log_{10}$-scale) in computing the double-layer
    potential using PTR$_{256}$ for the boundary
    $r(t) = {1.55 + 0.4 \cos 5 t }$ for the Dirichlet data,
    $f(\mathbf{y}) = \frac{1}{2\pi} \log | \mathbf{y} - \mathbf{x}_{0}
    |$
    with $\mathbf{x}_{0} = {(1.85, 1.65)}$. We evaluate the solution
    inside the domain on a regular grid.  [Right] Plot of the absolute
    error ($\log_{10}$-scale) in computing the
    double-layer potential using asymptotic
   PTR$_{\text{256}}$ given in
    \eqref{eq:DLP-asymptotic} for the same problem. The asymptotic
    method is used in a boundary layer determined by a threshold on
    the error from evaluating \eqref{kernel_identities}.}
  \label{fig:10}
\end{figure}

\begin{figure}[ht!]
  \centering
  \includegraphics[width=1\linewidth]{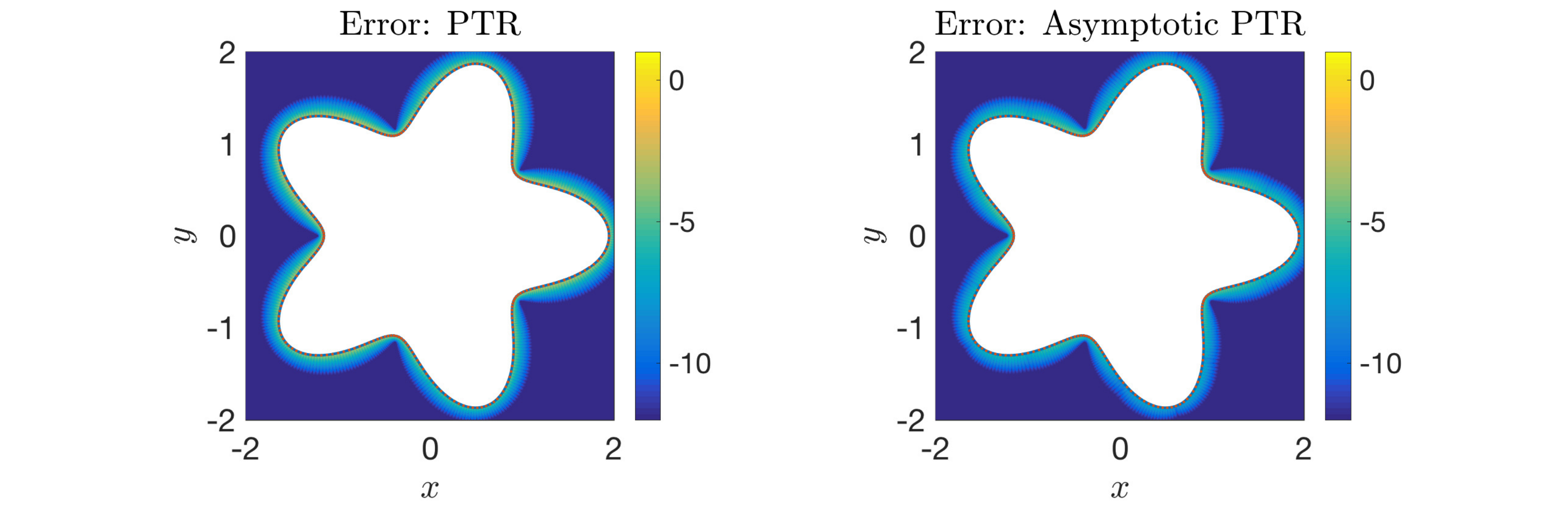}
  \caption{[Left] Plot of the absolute error
  ($\log_{10}$-scale) in computing the single-layer
    potential using PTR$_{256}$ for the boundary
    $r(t) = {1.55 + 0.4 \cos 5 t }$ for the Neumann data,
    $g(\mathbf{y}) = \frac{\partial v}{\partial \mathbf{n}}$ with
    $v(\mathbf{x}) = ( \mathbf{x} - \mathbf{x}_{0} ) / | \mathbf{x} -
    \mathbf{x}_{0} |^{2}$,
    $\mathbf{x}_{0} = {(0.1, 0.4)}$.  [Right] Plot of $\log_{10}$ of
    the absolute error ($\log_{10}$-scale) in
    computing the single-layer potential using the asymptotic
    PTR$_{256}$ given in \eqref{eq:4.15} for the same
    problem. The asymptotic method is used in a boundary layer
    determined by a threshold on the error from evaluating
    \eqref{kernel_identities}.}
      \label{fig:12}
\end{figure}

We present results of this generalized method for evaluation of the
double-layer potential and single-layer potential for the same
problems presented in Section \ref{ssec:DLPresults} and
\ref{ssec:SLPresults}, respectively.  We use a threshold of
$1 \times 10^{-8}$, as described above, to determine when the
evaluation point is inside the boundary layer and asymptotic
PTR$_{N}$ method is to be used.  In
{Fig.}~\ref{fig:10} we show the error in computing the double-layer
potential using PTR$_{256}$ and asymptotic
PTR$_{256}$ when solving on a Cartesian grid with
meshsize $h = 0.005$ within the boundary curve.  Similarly, we present
the evaluation of the single-layer potential in
{Fig.}~\ref{fig:12}. Here, we are computing on a Cartesian grid with
meshsize $h = 0.005$ exterior to the boundary curve. These results
show, similar to the results while considering a body-fitted grid,
that the error made by asymptotic PTR$_{N}$ is
several orders of magnitude smaller than those made by
PTR$_{N}$. However, there are more variations in
these errors because $\mathbf{y}^\ast$ does not always coincide with a
quadrature point.  We choose to use 256 quadrature points here as this
is what is actually needed to solve the boundary integral equations
for the densities such that $\mu(t) | \mathbf{y}'(t) |$ and
$\varphi(t) | \mathbf{y}'(t)|$ are sufficiently resolved. We were able
to use less points for the body-fitted grid as this restriction is not
as strict when $\mathbf{y}^\ast$ is a quadrature point.

\section{Conclusions}
\label{sec:conclusions}

We have presented a new method to address the close evaluation
problem.  When solving boundary value problems using boundary
integrals equation methods, the solution is evaluated at desired
points within the domain by numerically evaluating layer potentials.
Using the same quadrature rule that is used to solve the integral
equation for this evaluation achieves high order accuracy everywhere
in the domain, except close to the boundary where an $O(1)$ error is
incurred.  The new method developed here takes advantage of the
knowledge of the sharply peaked kernel of layer potentials close to
the boundary to reduce this error by several orders of magnitude. We
have used asymptotic methods to analyze the kernels of the double- and
single-layer potentials to relieve the numerical method from having to
integrate over this sharply peaked kernel. The resulting method is
straightforward to implement. We have presented results for both the
interior Dirichlet problem and exterior Neumann problem for Laplace's
equation and show a reduction in error of four to five orders of
magnitude in the solution evaluation close to the
boundary. Furthermore, we have presented how to generalize this method
to solve within the whole domain, including how to determine when to
use the new asymptotic method.

This asymptotic method has been recently applied to acoustic
scattering problems by sound-soft obstacles~\cite{carvalho2017}. For
those problems, the sharp peaks in the kernels for the single- and
double-layer potentials, both which are needed, have the same
character as those for Laplace's equation discussed here. Hence, only
small modifications are needed to apply these methods to wave
propagation problems.  We are currently extending this asymptotic
method to three-dimensional problems.  Furthermore, we are working on
different applications which include extending the method presented
here to the Stokes equations and studying surface plasmons.

\section*{Acknowledgments}

The authors thank Fran\c{c}ois Blanchette and Boaz Ilan for their
thoughtful discussions leading up to this work.  S. Khatri was
supported in part by the National Science
Foundation (PHY-1505061). A. D. Kim acknowledges support by Air
  Force Office of Scientific Research (FA9550-17-1-0238) and the
  National Science Foundation.

\end{document}